\documentclass[sn-mathphys]{sn-jnl}% Math and Physical Sciences Reference Style
\usepackage{mathrsfs}
\usepackage{amsmath}
\usepackage{amsfonts}
\usepackage{amsthm}
\usepackage{color}
\usepackage{xcolor}
\usepackage{graphicx}
\usepackage{tabularx}
\usepackage{makecell}
\usepackage{multirow}

\usepackage{bm}

\usepackage{caption}
\usepackage{subcaption}
\usepackage{float}
\usepackage{url}
\usepackage{multicol}
\usepackage{txfonts}
\usepackage{setspace}
\usepackage{arydshln}
\usepackage{enumitem}
\usepackage{lipsum}% for dummy text
\usepackage{siunitx,booktabs}
\usepackage{nth}
\usepackage{accents} % Accents: circle

\usepackage{listings}
\usepackage{xcolor}

\definecolor{codegreen}{rgb}{0,0.6,0}
\definecolor{codegray}{rgb}{0.5,0.5,0.5}
\definecolor{codepurple}{rgb}{0.58,0,0.82}
\definecolor{backcolour}{rgb}{0.95,0.95,0.92}

\lstdefinestyle{python-style}{
  backgroundcolor=\color{backcolour}, 
  commentstyle= \color{codegreen},
  keywordstyle=\color{black},
  numberstyle=\tiny\color{codegray},
  stringstyle=\color{codepurple},
  basicstyle=\ttfamily\footnotesize,
  breakatwhitespace=false,     
  breaklines=true,                 
  captionpos=b,                    
  keepspaces=true,                 
  showspaces=false,                
  showstringspaces=false,
  showtabs=false,                  
  tabsize=2,
  framexleftmargin=0.2cm,
  xleftmargin=0.5cm,xrightmargin=0.5cm
}

\labelformat{enumi}{#1}
\labelformat{enumii}{#1}
\labelformat{enumiii}{#1}
\labelformat{enumiv}{#1}

\jyear{2021}%

\theoremstyle{thmstyleone}%
\theoremstyle{thmstyletwo}%
\newtheorem{remark}{Remark}%

\theoremstyle{thmstylethree}%

\graphicspath{{figures/}}

\usepackage{listings}
\definecolor{light-gray}{gray}{0.95}
\begin{document}

\title[Open-source shape optimization of isogeometric shells]{Open-source shape optimization for isogeometric shells using FEniCS and OpenMDAO}

%%=============================================================%%
%% Prefix	-> \pfx{Dr}
%% GivenName	-> \fnm{Joergen W.}
%% Particle	-> \spfx{van der} -> surname prefix
%% FamilyName	-> \sur{Ploeg}
%% Suffix	-> \sfx{IV}
%% NatureName	-> \tanm{Poet Laureate} -> Title after name
%% Degrees	-> \dgr{MSc, PhD}
%% \author*[1,2]{\pfx{Dr} \fnm{Joergen W.} \spfx{van der} \sur{Ploeg} \sfx{IV} \tanm{Poet Laureate} 
%%                 \dgr{MSc, PhD}}\email{iauthor@gmail.com}
%%=============================================================%%

\author[1]{\mbox{\fnm{Han} \sur{Zhao}}}%\email{iauthor@gmail.com}
\author[1]{\mbox{\fnm{John T.} \sur{Hwang}}}
\author*[2,1]{\mbox{\fnm{Jiun-Shyan} \sur{Chen}}}\email{jsc137@ucsd.edu}

\affil[1]{\orgdiv{Department of Mechanical and Aerospace Engineering}, \orgname{University of California San Diego}, \orgaddress{\street{9500 Gilman Drive}, \city{La Jolla}, \state{CA} \postcode{92093}, \country{USA}}}

\affil[2]{\orgdiv{Department of Structural Engineering}, \orgname{University of California San Diego}, \orgaddress{\street{9500 Gilman Drive}, \city{La Jolla}, \state{CA} \postcode{92093}, \country{USA}}}

%%==================================%%
%% sample for unstructured abstract %%
%%==================================%%

\abstract{
We present an open-source Python framework for the shape optimization of complex shell structures using isogeometric analysis (IGA). IGA seamlessly integrates computer-aided design (CAD) and analysis models by employing non-uniform rational B-splines (NURBS) as basis functions, enabling the natural implementation of the Kirchhoff--Love shell model due to their higher order of continuity. We leverage the recently developed FEniCS-based analysis framework, PENGoLINS, for the direct structural analysis of shell structures consisting of a collection of NURBS patches through a penalty-based formulation. This contribution introduces the open-source implementation of gradient-based shape optimization for isogeometric Kirchhoff--Love shells with a modular architecture. Complex shell structures with non-matching intersections are handled using a free-form deformation (FFD) approach and a moving intersections formulation. The symbolic differentiation and code generation capabilities in FEniCS are utilized to compute the analytical derivatives. By integrating FEniCS with OpenMDAO, we build modular components that facilitate gradient-based shape optimization of shell structures. The modular architecture in this work supports future extensions and integration with other disciplines and solvers, making it highly customizable and suitable for a wide range of applications. We validate the design-analysis-optimization workflow through several benchmark problems and demonstrate its application to aircraft wing design optimization. The framework is implemented in a Python library named GOLDFISH (Gradient-based Optimization and Large-scale Design Framework for Isogeometric SHells) and the source code will be maintained at \url{https://github.com/hanzhao2020/GOLDFISH}.}


%%\pacs[JEL Classification]{D8, H51}

%%\pacs[MSC Classification]{35A01, 65L10, 65L12, 65L20, 65L70}

\maketitle

% % Statements and highlights
% \section*{Statements and declarations}
% \subsection*{Conflict of interest}
% The authors have no relevant financial or non-financial interests to disclose.

% \section*{Article highlights}
% \begin{itemize}
%     \item Open-source implementation of immersogeometric fluid--structure interaction using FEniCS-based code generation.
%     \item Background fluid domain deforms, using a symbolic pullback to a reference domain.
%     \item Code generation enables prototyping of constitutive models for immersed structure.
% \end{itemize}

\section{Introduction} \label{sec:introduction}
Shape optimization of shell structures is crucial in engineering design to improve structural performance, efficiency, and material utilization \cite{bletzinger1993form}. Traditional shape optimization methods using the finite element method (FEM) often struggle with accurately representing complex geometries and generating high-quality meshes \cite{Hardwick2005}. Inaccurate geometry representation introduces errors in the analysis and optimization, while poor mesh quality leads to numerical issues. To address these difficulties, shape optimization using IGA \cite{Hughes05a, CoHuBa09} has emerged as a powerful method. IGA offers seamless integration between design geometries and analysis models by adopting the same mathematical description. By utilizing NURBS \cite{piegl2012nurbs} or other types of splines \cite{sederberg2003t, giannelli2012thb, thomas2022u} as basis functions, IGA allows direct analysis of CAD geometries with improved accuracy \cite{Evans2009}. This property makes IGA particularly promising for the analysis and optimization of shell structures, enabling updated design geometries to be used directly for structural analysis without the need for finite element (FE) mesh generation. The natural fulfillment of the $C^1$ continuity requirement in the Kirchhoff--Love shell model \cite{Kiendl2009, Kiendl2011} further enhances the applicability of IGA to thin-walled structures, which are widespread in structural design.

For real-world complex shell structures typically modeled by multiple NURBS surfaces, patch coupling at surface intersections is needed in analysis and optimization. Maintaining continuities of displacements and rotations at these intersections is essential. Various methods have been proposed to couple separate shell patches, including the bending strip method \cite{Bazilevs10c} for conforming discretizations. Nitsche-type formulations \cite{guo2015nitsche, Guo2018, guo2021isogeometric, Benzaken2021, wang2022isogeometric}, mortar methods \cite{BRIVADIS2015292, horger2019hybrid, hirschler2019dual, schuss2019multi}, and penalty-based methods \cite{herrema2019penalty,leonetti2020robust,proserpio2022penalty,zhao2022open, guarino2024interior} are developed for coupling of Kirchhoff--Love shell patches with non-conforming intersections. Furthermore, \cite{farahat2023isogeometric} introduced a strong coupling approach by approximating a collection of shell patches as an analysis-suitable $G^1$ multi-patch surface. Additionally, T-splines \cite{bazilevs2010isogeometric, casquero2020seamless} and G-splines \cite{wen2023isogeometric} have been utilized for the structural analysis of complex isogeometric shell structures.

Shape optimization for non-matching shell structures poses additional technical challenges, as the intersections must be managed during the optimization process. Without proper handling, initial intersecting shell patches may become separated or self-penetrated, leading to unrealistic structures. The spline composition method \cite{hirschler2019embedded, hirschler2019isogeometric, bouclier2022iga, hao2023isogeometric} was proposed for maintaining the surface intersection during shape optimization but requires identifying master surfaces and extruding a 3D solid from them. An FFD \cite{Sederberg1986}-based shape optimization approach \cite{zhao2024automated} ensures the connectivity of shell patches at non-matching intersections but may suffer from significant element distortion with large intersection movements. Recently, \cite{zhao2024shape} proposed an optimization method to overcome the element distortion issue, allowing relative movement between intersecting shell patches while retaining intersections.

Code transparency in the field of design optimization has been gaining increasing interest. OpenMDAO \cite{gray2019openmdao}, an open-source framework for multidisciplinary design optimization (MDO), uses a modular architecture that allows users to create custom models for different disciplines and integrate them with various optimization algorithms. Successful applications of OpenMDAO span aerospace engineering \cite{Jasa2018a,Jasa2020,Adler2022b}, wind energy \cite{Herrema2019a}, robotics \cite{yan2019,lin2022generalized}, and topology optimization \cite{chung2019topology, yan2022topology, jauregui2023}. The Python library CSDL \cite{gandarillas2024graph} addresses large-scale MDO problems using a graph representation to automatically generate adjoint sensitivities. The CSDL-based Python library FEMO \cite{xiang2024}, coupled with the FEniCS project \cite{Logg:2010, Logg2012, Alnaes2014}, was developed to solve partial differential equation (PDE)-constrained optimization problems, significantly reducing the coding effort for PDE components such as structural mechanics, fluid mechanics, and heat transfer, etc. Despite these advancements, an open-source design optimization framework using IGA has been lacking. This contribution fills that gap by developing an open-source Python library for shape optimization of complex shell structures using IGA, enabling researchers to explore the benefits of IGA in shape optimization problems and advance structural design optimization.

In the proposed framework GOLDFISH, we use OpenMDAO as the optimization toolkit, with plans to incorporate CSDL in future work to make use of the graph-based paradigm. For the IGA solver for structural analysis of non-matching shell structures, we employ the open-source package PENGoLINS \cite{zhao2022open}. PENGoLINS, a Python framework based on tIGAr \cite{Kamensky2019}, uses extraction techniques \cite{BSEH11, SBVSH11, Schillinger2016, fromm2023interpolation, fromm2024interpolation} to construct spline basis functions from Lagrange basis functions in the FEM solver FEniCS. tIGAr has been successfully applied in various fields \cite{Bazilevs2019, Eikelder2019, Zhang2019, yang2020determination}. Additionally, an open-source fluid--structure interaction framework \cite{Kamensky2021, neighbor2023leveraging} is developed based on tIGAr and shows a good application for prosthetic heart valve simulation with isogeometric leaflets. PENGoLINS employs a penalty-based formulation \cite{zhao2022open} to couple the non-matching isogeometric Kirchhoff--Love shells, which is automated using the code generation technology in FEniCS. The current code framework incorporates IGA for complex Kirchhoff--Love shells with a penalty formulation \cite{herrema2019penalty}, and utilizes the FFD-based and moving intersections approaches \cite{zhao2024automated, zhao2024shape} to handle patch intersections. The modular architecture of OpenMDAO ensures that GOLDFISH can readily couple with other disciplines and extend to more practical problems.

The remainder of this paper is outlined as follows. Section \ref{sec:formulation} provides an overview of shape optimization approaches for non-matching isogeometric shell structures. Section \ref{sec:design-goldfish} discusses technical details and framework structures of GOLDFISH. A series of benchmarks with code implementation are presented in Section \ref{sec:numerical-examples} to validate the framework, alongside demonstrations of aircraft wing applications. Section \ref{sec:conclusions} draws conclusions and discusses future work.

~\\

\section{Optimization formulations} \label{sec:formulation}
This section presents the shape optimization schemes for isogeometric shells employed in GOLDFISH. We first review the basic isogeometric Kirchhoff--Love shell formulations in Section \ref{subsec:isogeometric-KL-shells}, then introduce the overall shape optimization strategy for isogeometric shells in Section \ref{subsec:shape-opt-isogeometric-shells}. In Section \ref{subsec:shape-opt-non-matching-shells}, we discuss methods for handling patch intersections in complex shell structures.

\subsection{Kirchhoff--Love shell with isogeometric discretization} \label{subsec:isogeometric-KL-shells}
The basic derivation of the Kirchhoff--Love shell is presented in this Section to provide the foundation for the subsequent shape optimization formulations. The geometry and unknowns of the Kirchhoff--Love shell are discretized isogeometrically. The Kirchhoff--Love shell theory assumes negligible transverse shear strains and constant thickness during deformation, allowing the shell to be represented by its mid-surface geometry $\mathbf{X} \in \mathbb{R}^d$, with the deformation described by $\mathbf{u}\in \mathbb{R}^d$, where $d$ is the spatial dimension. For an input CAD geometry of the mid-surface parametrized with coordinates $\bm{\xi} = \{ \xi_\alpha \}_{\alpha=1}^{2}$ and discretized by B-spline or NURBS basis functions $\mathbf{N}(\bm{\xi}) = \begin{bmatrix} N_1(\bm{\xi}) & N_2(\bm{\xi}) & \ldots & N_n(\bm{\xi}) \end{bmatrix}$, where $n$ is the number of unknowns, the mid-surface geometry and displacement field are expressed as
\begin{align}
    \mathbf{X}(\bm{\xi}) = \sum\limits_{i=1}^{n} {N}_i(\bm{\xi}) {P}_i = \mathbf{N}^{\text{T}}(\bm{\xi}) \mathbf{P} \quad \text{ and } \quad \mathbf{u}(\bm{\xi}) = \sum\limits_{i=1}^{n} {N}_i(\bm{\xi}) {d}_i = \mathbf{N}^{\text{T}}(\bm{\xi}) \mathbf{d} \text{ ,} \label{eq:geom-disp-discretization}
\end{align}
where $\mathbf{P}$ and $\mathbf{d}$ are the geometric and displacement control points, respectively. The deformed state of the shell is defined as
\begin{align}
    \mathbf{x}(\bm{\xi}) = \mathbf{X}(\bm{\xi}) + \mathbf{u} (\bm{\xi}) = \mathbf{N}^{\text{T}}(\bm{\xi}) (\mathbf{P} + \mathbf{d}) \text{ .}
\end{align}
The local covariant basis vectors in the reference and deformed configurations are formulated as
\begin{align}
    \mathbf{A}_{\alpha} = \mathbf{X},_{\xi_\alpha} = \mathbf{N},_{\xi_\alpha}^{\text{T}} \mathbf{P} \quad \text{ and } \quad \mathbf{a}_{\alpha} = \mathbf{x},_{\xi_\alpha} = \mathbf{N},_{\xi_\alpha}^{\text{T}} (\mathbf{P} + \mathbf{d}) \text{ ,}
\end{align}
and the unit normal vectors are given by
\begin{align}
    \mathbf{A}_3 = \frac{\mathbf{A}_1\times \mathbf{A}_2}{\Vert \mathbf{A}_1\times \mathbf{A}_2 \Vert} \quad \text{ and } \quad \mathbf{a}_3 = \frac{\mathbf{a}_1\times \mathbf{a}_2}{\Vert \mathbf{a}_1\times \mathbf{a}_2 \Vert} \text{ ,}
\end{align}
where $\Vert \cdot \Vert$ is the Euclidean norm. The metric and curvature coefficients in the reference configuration are given by
\begin{align}
    A_{\alpha \beta} = \mathbf{A}_{\alpha} \cdot \mathbf{A}_{\beta} \quad \text{ and } \quad B_{\alpha \beta} = \mathbf{A}_{\alpha},_{\xi_{\beta}} \cdot \mathbf{A}_3 =  - \mathbf{A}_{\alpha} \cdot \mathbf{A}_3,_{\xi_{\beta}} \text{ .}
\end{align}
Analogously, the associated coefficients in the deformed configuration are
\begin{align}
    a_{\alpha \beta} = \mathbf{a}_{\alpha} \cdot \mathbf{a}_{\beta} \quad \text{ and } \quad b_{\alpha \beta} = \mathbf{a}_{\alpha},_{\xi_{\beta}} \cdot \mathbf{a}_3 =  - \mathbf{a}_{\alpha} \cdot \mathbf{a}_3,_{\xi_{\beta}} \text{ .}
\end{align}
The coefficients of the membrane strain tensor and curvature change tensor are defined as
\begin{align}
    \varepsilon_{\alpha \beta} = \frac{1}{2}(a_{\alpha \beta} - A_{\alpha \beta}) \quad \text{ and } \quad \kappa_{\alpha \beta} = B_{\alpha \beta} - b_{\alpha \beta} \text{ .}
\end{align}
By organizing the membrane strain tensor and curvature change tensor using Voigt notation as $\bm{\varepsilon}$ and $\bm{\kappa}$, the normal forces and bending moments of the shell, modeled with the St. Venant--Kirchhoff material model, read as
\begin{align}
    \mathbf{n} = t\, \mathbf{C} : \bm{\varepsilon} \quad \text{ and } \quad \mathbf{m} = \frac{t^3}{12}\, \mathbf{C} : \bm{\kappa} \text{ ,} 
\end{align}
where $t$ is the shell thickness and $\mathbf{C}$ is the material tensor.
% They can be obtained by discretizing the partial derivative of the total energy with respect to shell displacement
% \begin{align}
%     \mathbf{R}_{\text{S}} =  \partial_{\mathbf{d}} W_{\text{S}} \text{ .} \label{eq:KL-shell-residual}
% \end{align}
The total energy of the Kirchhoff--Love shell theory is stated as
\begin{align}
    W_{\text{S}} = W^{\text{int}} - W^{\text{ext}} = \int_S \mathbf{n} : \bm{\varepsilon} + \mathbf{m} : \bm{\kappa}  \, \mathrm{d}S -  \int_S \mathbf{f} \cdot \mathbf{u} \, \mathrm{d}S \text{ ,} \label{eq:total-energy-KL-shell}
\end{align}
where $\mathbf{f}$ stands for the external force applied on $\mathbf{S}$. For a detailed derivation of the Kirchhoff--Love shell model, readers are referred to \cite[Section 3]{Kiendl2011}. The curvature change tensor and associated bending moments involve the second-order derivative of the displacement, which requires the basis functions in the solution space to be $C^1$ continuous across element boundaries. This requirement is naturally satisfied by the B-spline and NURBS basis functions, which demonstrate excellent results \cite{Kiendl2009, Kiendl2015, Kamensky2015}. 

\begin{remark}
Although shear locking is avoided in the Kirchhoff--Love shell theory, membrane locking remains a crucial numerical challenge for curved thin shells. While in the current framework, numerous numerical studies have been conducted in \cite{zhao2022open} and \cite{zhao2024automated} on several benchmark problems to evaluate the accuracy and convergence of the analysis code, implementing methods to eliminate membrane locking in isogeometric Kirchhoff--Love shells into this open-source framework can be considered. Several classes of methods are applicable. The first class belongs to the strain projection technique, such as the assumed strain method \cite{casquero2022removing, casquero2023overcoming, mathews2024computationally}, B-bar type formulation \cite{bouclier2013efficient, greco2018reconstructed}, and the discrete strain gap technique \cite{koschnick2005discrete, echter2013hierarchic}. These methods can be implemented directly under the present code structure. The hybrid discretization method in \cite{sauer2024simple} with different orders of continuity for membrane and bending strains can also be implemented in the present framework with ease. The other class of methods involving multi-field variational principles, such as \cite{bieber2018variational} based on the Hellinger--Reissner variational principle, require changes of the present code structure due to additional unknowns in the formulation.
\end{remark}

% Equation \eqref{eq:KL-shell-residual} represents a system of nonlinear equations, which can be linearized and solved by Newton--Raphson method for the displacement increments
% \begin{align}
%     \partial_{\mathbf{d}} \mathbf{R}_{\text{S}} \, \Delta \mathbf{d} = -\mathbf{R}_{\text{S}} \text{ .} \label{eq:KL-shell-linearization}
% \end{align}
% The discretization of $\partial_{\mathbf{d}} \mathbf{R}_{\text{S}}$ represents the stiffness matrix $\mathbf{K}_{\text{S}} = \partial_{\mathbf{d}} \mathbf{R}_{\text{S}}$ of the shell structure.

\subsection{Shape optimization for isogeometric shells} \label{subsec:shape-opt-isogeometric-shells}
The shape optimization for shell structures can be formulated as the following form
\begin{equation}
\begin{aligned}
    \underset{\mathbf{P}}{\text{minimize }} &f(\mathbf{P}) \\
    \text{subject to } &\mathbf{g}(\mathbf{P}) \leq \mathbf{0}\\
    &\mathbf{h}(\mathbf{P}) = \mathbf{0} \text{ ,} \label{eq:minmize-single-patch}
\end{aligned}
\end{equation}
where $\mathbf{P}$ are the design variables, $f$ is the objective function, and $\mathbf{g}$ and $\mathbf{h}$ represent the inequality and equality constraints, respectively. In shape optimization problems using IGA, the design variables $\mathbf{P}$ are the coordinates of the control points that define the shell geometry using B-spline or NURBS basis functions. To enable gradient-based design optimization, it is necessary to derive the total derivative of the objective function with respect to the design variables \cite{kiendl2014isogeometric, hirschler2021new}. For shape optimization, the model outputs are typically affected by the displacement of the structure, which depends on the design variables and is referred to as the state variable. Therefore, we express the objective function as $f(\mathbf{P}, \mathbf{d}(\mathbf{P}))$, and the total derivative is given by
\begin{align}
    \mathrm{d}_{\mathbf{P}} f = \partial_{\mathbf{P}}f + \left( \partial_{\mathbf{d}}f \right)^{\mathrm{T}}  \mathrm{d}_{\mathbf{P}}\mathbf{d} \text{ ,} \label{eq:derivative-df-dp}
\end{align}
where partial derivatives $\partial_{\mathbf{P}}f$ and $\partial_{\mathbf{d}}f$ can be readily computed. However, to obtain the total derivative in \eqref{eq:derivative-df-dp}, we need to derive the total derivative of the structural displacement with respect to shell control points $\mathrm{d}_{\mathbf{P}}\mathbf{d}$. In the context of shape optimization for shell structures, the relation between the geometry control points $\mathbf{P}$ and the displacement $\mathbf{d}$ is typically governed by a system of equations derived from discretized PDEs
\begin{align}
    \mathbf{R}_{\text{S}}(\mathbf{P}, \mathbf{d}) = \mathbf{0} \text{ ,}
\end{align}
where $\mathbf{R}_{\text{S}}$ represents the residual equations. They can be obtained by discretizing the partial derivative of the total energy in \eqref{eq:total-energy-KL-shell} with respect to shell displacement
\begin{align}
    \mathbf{R}_{\text{S}} = \partial_{\mathbf{d}} W_{\text{S}} \text{ .} \label{eq:KL-shell-residual}
\end{align}
Equation \eqref{eq:KL-shell-residual} represents a system of nonlinear equations, which can be linearized and solved by the Newton--Raphson method to compute the displacement increments
\begin{align}
    \partial_{\mathbf{d}} \mathbf{R}_{\text{S}} \, \Delta \mathbf{d} = -\mathbf{R}_{\text{S}} \text{ .} \label{eq:KL-shell-linearization}
\end{align}
The discretization of $\partial_{\mathbf{d}} \mathbf{R}_{\text{S}}$ represents the stiffness matrix $\mathbf{K}_{\text{S}} = \partial_{\mathbf{d}} \mathbf{R}_{\text{S}}$ of the shell structure.

With the adjoint method, the total derivative $\mathrm{d}_{\mathbf{P}}\mathbf{d}$ in \eqref{eq:derivative-df-dp} can be obtained by taking the total derivative of $\mathbf{R}_{\text{S}}$ with respect to $\mathbf{P}$
\begin{align}
    \mathrm{d}_{\mathbf{P}}\mathbf{R}_{\text{S}} = \partial_{\mathbf{P}}\mathbf{R}_{\text{S}} + \partial_{\mathbf{d}}\mathbf{R}_{\text{S}}\, \mathrm{d}_{\mathbf{P}}\mathbf{d} = \mathbf{0} \text{ ,} \label{eq:derivative-dr-dp}
\end{align}
and substituting $\partial_{\mathbf{d}}\mathbf{R}_{\text{S}}$ by $\mathbf{K}_{\text{S}}$
\begin{align}
    \mathrm{d}_{\mathbf{P}}\mathbf{d} = - \mathbf{K}_{\text{S}}^{-1}\,\partial_{\mathbf{P}}\mathbf{R}_{\text{S}} \text{ .} \label{eq:derivative-du-dp}
\end{align}
The partial derivative $\partial_{\mathbf{P}}\mathbf{R}_{\text{S}}$ can be computed from the residual vector in a similar way to the stiffness matrix.

Substituting \eqref{eq:derivative-du-dp} into \eqref{eq:derivative-df-dp}, the total derivative for the shape optimization problem of an isogeometric Kirchhoff--Love shell can be expressed as
\begin{align}
    \mathrm{d}_{\mathbf{P}} f = \partial_{\mathbf{P}}f - \left( \partial_{\mathbf{d}}f \right)^{\mathrm{T}}  \mathbf{K}_{\text{S}}^{-1}\,\partial_{\mathbf{P}}\mathbf{R}_{\text{S}}  \text{ .} \label{eq:derivative-df-dp-full}
\end{align}
Using the total derivative, we can apply gradient-based optimization algorithms \cite{ahmadianfar2020gradient} to determine the optimal shape for a given baseline design of a shell structure.

\subsection{Treatment of surface intersections} \label{subsec:shape-opt-non-matching-shells}
\begin{figure}[!htb]\centering
    \includegraphics[height=1.5in]{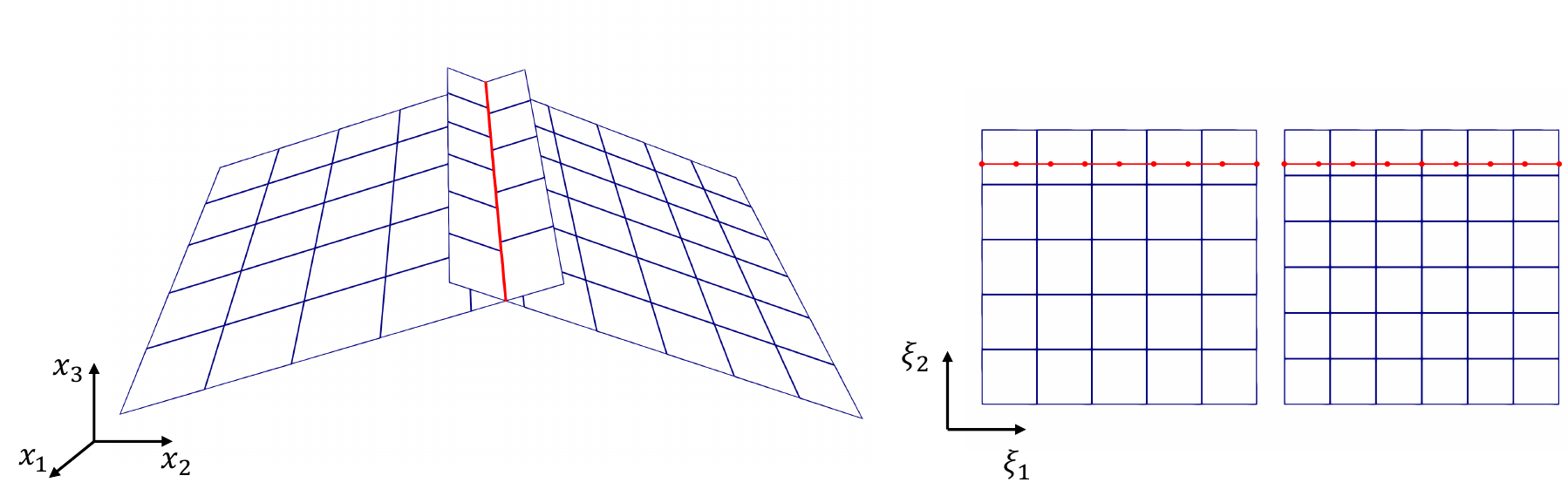}
    \caption{An illustrative example of two shell patches with one intersection, where shell patches are described by NURBS surfaces. A topologically 1D, geometrically 2D quadrature mesh is generated in the parametric space to enforce displacement and rotational continuity through a penalty formation.}
    \label{fig:shell-coupling}
\end{figure}

CAD geometries are typically modeled by multiple NURBS patches for real-world shell structures. The multi-patch NURBS-based geometries exhibit intersections between patches, where displacement and angular compatibility needs to be maintained during structural analysis to glue the patches together. An illustrative example of two shell patches with one intersection is depicted in Figure \ref{fig:shell-coupling}. Additionally, it is crucial to preserve compatibility at patch intersections during shape updates in the optimization loop so that the optimized shell structures will not improperly separate or result in unrealistic structures. The following sections discuss the methods used to handle shape optimization for shell structures consisting of multiple spline patches.

\subsubsection{Compatibility conservation of surface intersections using FFD} \label{subsubsec:shape-opt-ffd}
\begin{figure}[!htb]\centering
    \includegraphics[height=1.6in]{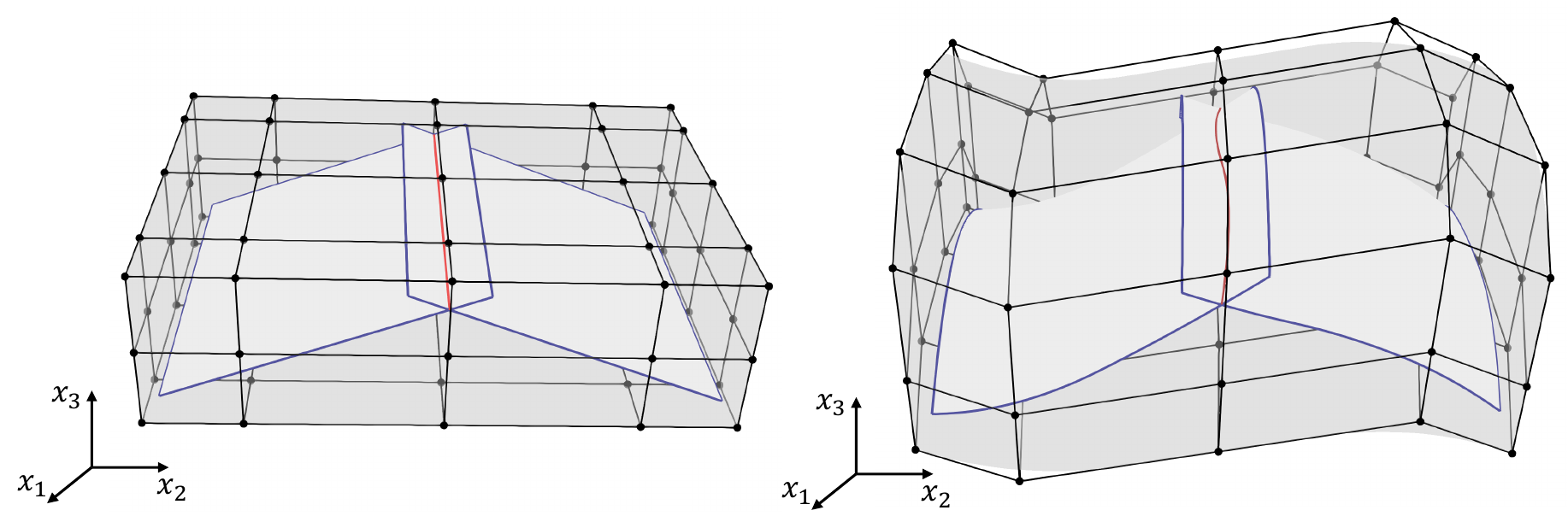}
    \caption{Illustration of the FFD-based shape optimization approach. The intersecting shell patches are embedded in a trivariate B-spline FFD block, whose control points serve as design variables. Modifications to the FFD block shape alter the Lagrange nodal points of the shell patches, while preserving patch intersections. The updated NURBS control points are then obtained from the Lagrange nodal points using Lagrange extraction matrices.}
    \label{fig:shell-shopt-ffd}
\end{figure}

In the proposed code framework, we implement the FFD-based shape optimization formulation \cite{zhao2024automated}, combined with the Lagrange extraction technique \cite{Schillinger2016}, to preserve intersections among separately modeled spline surfaces within a single CAD geometry. In this approach, the entire CAD geometry of a shell structure is embedded in a trivariate B-spline block. The control points of the FFD block are related to the Lagrange nodal points of the shell geometry interpolated with Lagrange polynomials, preserving patch intersections inside the FFD block while modifying the shell geometry. The Lagrange nodal points of the shell structure are updated through the shape changes in the 3D FFD block during shape optimization. Subsequently, a pseudo-inverse system is solved to approximate the control points of the NURBS geometry from the updated Lagrange nodal points using the Lagrange extraction matrices. This approach is illustrated in Figure \ref{fig:shell-shopt-ffd}, where the shape of shell patches follows the shape change of the FFD block. As the shape changes continuously inside the FFD block, the embedded surfaces maintain their intersections without relative movement. Therefore, the control points of the 3D block are considered design variables.

Using the Lagrange extraction method, the NURBS basis functions of spline patch $k$ can be expressed in terms of Lagrange polynomials $\mathbf{N}^{k}_{\text{L}}$,
\begin{equation}
    \mathbf{N}^{k} = {\mathbf{M}^{k}}^{\text{T}} \, \mathbf{N}^{k}_{\text{L}} \text{ ,}
\end{equation}
where $\mathbf{M}^{k}$ denotes the extraction matrix. Each entry of $\mathbf{M}^{k}$ is the evaluation of the NURBS basis functions $\mathbf{N}^{k}$ at the parametric location of a Lagrange nodal point $\bm{\xi}^{k}$,
\begin{align}
    {M}^{k}_{ij} = {N}_j^{k} (\bm{\xi}^{k}_i) \label{eq:lagrange-extraction-matrix-generation} \text{ .}
\end{align}
The relationship between shell patches' NURBS control points $\mathbf{P}$ and Lagrange nodal points $\mathbf{P}_{\text{L}}$ for the $k$-th shell patch is defined by the Lagrange extraction matrix
\begin{align}
    \mathbf{M}^{k} \mathbf{P}^{k} = {\mathbf{P}^{k}_{\text{L}}} \text{ ,} \label{eq:lagrange-extraction}
\end{align}
It is noted that the extraction operator $\mathbf{M}^{k}$ is a nonsquare matrix and a pseudo-inverse is necessary to obtain the total derivative of the $k$-th spline patch,
\begin{align}
    \mathrm{d}_{\mathbf{P}^{k}_{\text{L}}} \mathbf{P}^{k} = \left({\mathbf{M}^{k}}^{\text{T}} \mathbf{M}^{k}\right)^{-1} {\mathbf{M}^{k}}^{\text{T}} \text{ .} \label{eq:derivative-dpk-dplk}
\end{align}
And the total derivative for the entire shell structure with $m$ spline patches is 
\begin{align}
    \mathrm{d}_{\mathbf{P}_{{\text{L}}}} \mathbf{P} = \text{diag}\left(\mathrm{d}_{\mathbf{P}^{1}_{\text{L}}} \mathbf{P}^{1}, \mathrm{d}_{\mathbf{P}^{2}_{\text{L}}} \mathbf{P}^{2}, \ldots, \mathrm{d}_{\mathbf{P}^{m}_{\text{L}}} \mathbf{P}^{m}\right) \text{ .} \label{eq:derivative-dp-dpl}
\end{align}

The total derivative between the control points of the FFD block and Lagrange nodal points of the $k$-th shell patch can be derived in a similar way
\begin{align}
    \mathbf{A}^{k} \mathbf{P}_{\text{FFD}} = \mathbf{P}^{k}_{\text{L}} \text{ ,} \label{eq:cp-ffd-cp-lagrange}
\end{align}
and the matrix $\mathbf{A}^{k}$ is constructed as
\begin{align}
    {A}_{ij}^{k} = N_{\text{FFD}j}(\mathbf{P}^{k}_{\text{L}i}) \text{ ,} \label{eq:derivative-dpl-dpffd-single-patch}
\end{align}
where $\mathbf{N}_{\text{FFD}}$ are the B-spline basis functions of the FFD block and the identity geometric mapping for the FFD block is assumed for simplicity. The associated total derivative $\mathrm{d}_{\mathbf{P}_{\text{FFD}}} \mathbf{P}_{\text{L}}$ is expressed as
\begin{align}
    \mathrm{d}_{\mathbf{P}_{\text{FFD}}} \mathbf{P}_{\text{L}} = \text{diag}\left( \mathbf{A}^{1}, \mathbf{A}^{2}, \ldots, \mathbf{A}^{m} \right) \text{ .} \label{eq:derivative-dpl-dpffd}
\end{align}

Analogous to \eqref{eq:derivative-df-dp-full}, the total derivative of the FFD-based shape optimization approach for a multi-patch shell structure is expressed as
\begin{align}
    \mathrm{d}_{\mathbf{P}_{\text{FFD}}} f = \left( \partial_{\mathbf{P}}f - \left( \partial_{\mathbf{d}}f \right)^{\mathrm{T}}  \mathbf{K}^{-1}\,\partial_{\mathbf{P}}\mathbf{R} \right) \mathrm{d}_{\mathbf{P}_{\text{L}}} \mathbf{P} \, \mathrm{d}_{\mathbf{P}_{\text{FFD}}} \mathbf{P}_{\text{L}} \text{ ,} \label{eq:derivative-df-dpffd-full}
\end{align}
where $\mathbf{P}^{\mathrm{T}} = \begin{bmatrix}
    {\mathbf{P}^{1}}^{\mathrm{T}} & {\mathbf{P}^{2}}^{\mathrm{T}} & \ldots & {\mathbf{P}^{m}}^{\mathrm{T}}
\end{bmatrix}$ is the B-spline or NURBS control points of all shell patches in the CAD geometry, and $\mathbf{P}_{\text{L}}^{\mathrm{T}} = \begin{bmatrix}
    {\mathbf{P}_{\text{L}}^{1}}^{\mathrm{T}} & {\mathbf{P}_{\text{L}}^{2}}^{\mathrm{T}} & \ldots & {\mathbf{P}_{\text{L}}^{m}}^{\mathrm{T}}
\end{bmatrix}$ is the Lagrange nodal points of the shell geometry. ${\mathbf{P}_{\text{FFD}}}$ denotes the control points of the B-spline FFD block. And $\mathbf{d}^{\mathrm{T}} = \begin{bmatrix}
    {\mathbf{d}^{1}}^{\mathrm{T}} & {\mathbf{d}^{2}}^{\mathrm{T}} & \ldots & {\mathbf{d}^{m}}^{\mathrm{T}}
\end{bmatrix}$ represents the displacements of all shell patches. Additionally, $\mathbf{R}(\mathbf{P}, \mathbf{d})$ is the residual vector of the multi-patch shell structure, and $\mathbf{K}$ denotes the associated stiffness matrix.

For the structural analysis of shell geometries consisting of a collection of NURBS surfaces, we employ a penalty coupling formulation to maintain displacement and rotational compatibility at the intersections. The residual vector and stiffness matrix of the multi-patch shell structure are expressed as 
\begin{align}
    \mathbf{R} = \partial_{\mathbf{d}} \left( \sum\limits^{m}_{k=1} W^{k}_{\text{S}} + \sum\limits^{r}_{l=1} W^{l}_{\text{pen}} \right) \quad \text{ and } \quad \mathbf{K} = \partial_{\mathbf{d}} \mathbf{R} \text{ ,} \label{eq:nonmatching-residual-stiffness}
\end{align}
where $W_{\text{pen}}$ is the penalty energy between two intersecting shell patches as illustrated in Figure \ref{fig:shell-coupling} and $r$ is the number of intersections. The penalty energy of the $l$-th intersection between shell patches with indices $l_1$ and $l_2$ is expressed as
\begin{equation}
\begin{aligned}
    %W^{l}_\text{pen} = \frac{1}{2}\int_{\mathcal{L}}\alpha_\text{d}\left\Vert\mathbf{d}^{l_1}-\mathbf{d}^{l_2}\right\Vert^2 + \alpha_\text{r}\left(\left(\mathbf{a}_3^{l_1}\cdot\mathbf{a}_3^{l_2} - \mathbf{A}_3^{l_1}\cdot{\mathbf{A}}_3^{l_2}\right)^2 + \left(\mathbf{a}_n^{l_1}\cdot\mathbf{a}_3^{l_2} - \mathbf{A}_n^{l_1}\cdot{\mathbf{A}}_3^{l_2}\right)^2\right)\,d\mathcal{L}\text{ ,} \label{eq:penalty-energy}
    W^{l}_{\text{pen}} & = \int_{\mathcal{L}} \alpha_d \Vert \mathbf{u}^{l_1} - \mathbf{u}^{l_2} \Vert^2 +  \alpha_r \left( (\sin \phi - \sin \phi_0 )^2 + (\cos \phi - \cos \phi_0 )^2 \right) \, \mathrm{d}\mathcal{L} \text{ ,} \label{eq:penalty-energy}
\end{aligned}
\end{equation}
where $\mathcal{L}$ represents the intersection curve with associated displacements $\mathbf{u}^{l_1}$ and $\mathbf{u}^{l_2}$, and $\phi_0$ and $\phi$ denote the angles between the two shell patches before and after deformation, respectively. The first term in \eqref{eq:penalty-energy} ensures that the two intersecting shells have the same displacement along the intersection, while the second term maintains the angle between the two shell patches during deformation. The penalty parameters $\alpha_d$ and $\alpha_r$ can be determined based on the material properties of the problem and shell discretizations. Further details on the penalty formulation can be found in \cite[Section 2]{herrema2019penalty}, and the open-source implementation using FEniCS is presented in \cite{zhao2022open}.

Substituting \eqref{eq:derivative-dp-dpl} and \eqref{eq:derivative-dpl-dpffd} into \eqref{eq:derivative-df-dpffd-full}, we can obtain the total derivative of the FFD-based approach for non-matching shell structures to perform shape optimization. The parametric locations of surface intersections are maintained by embedding the shell patches into the FFD block and there is no relative movement between shell patches.

\subsubsection{Shape updates with moving intersections}  \label{subsubsec:shape-opt-moving-int}
\begin{figure}[!htb]\centering
    \includegraphics[height=1.5in]{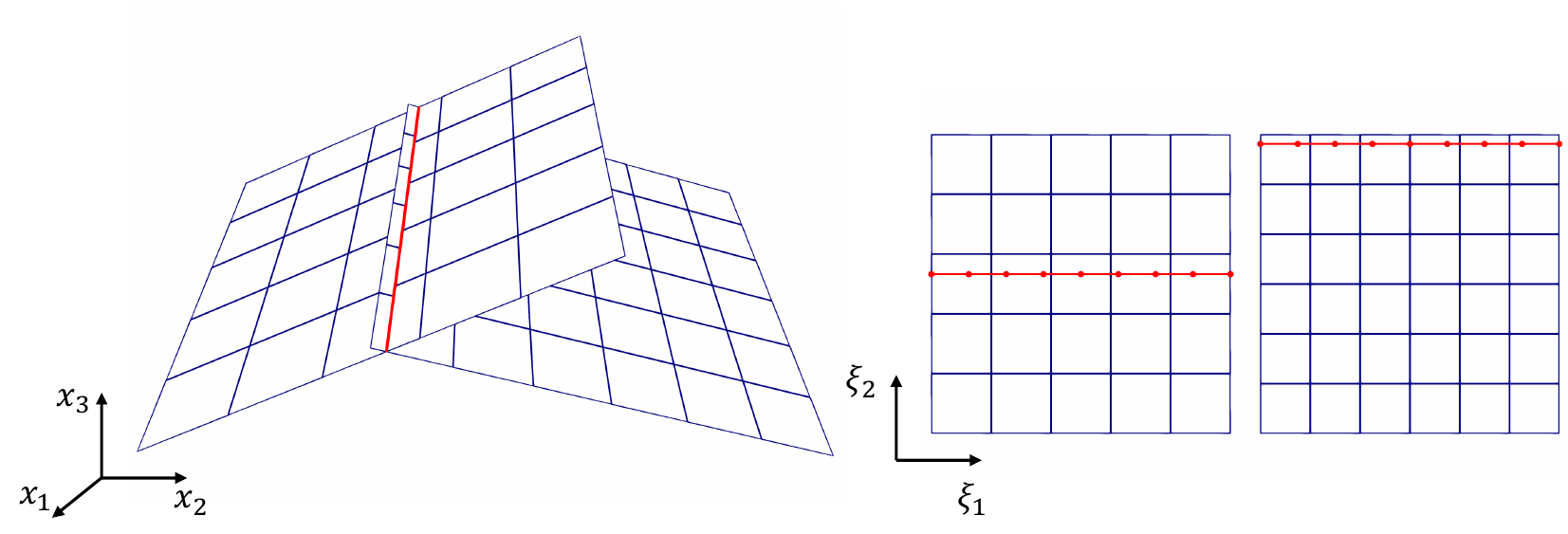}
    \caption{Updated design of the intersecting patches from Figure \ref{fig:shell-coupling}. The shape of the two shell patches is allowed to change relative to each other, and parametric coordinates of the intersection are determined by solving the implicit equation \eqref{eq:intersection-residual}.}
    \label{fig:shell-shopt-mint}
\end{figure}

While the FFD-based approach is effective in preserving surface intersections in shell structural optimization, it may significantly deteriorate the quality of shell elements if shape optimization involves substantial movement of patch intersections. To address this issue, \cite{zhao2024shape} proposed an optimization approach for shell structures with moving intersections. In this approach, shell patches are allowed to move relative to other intersecting patches in the parametric coordinates during the shape update process without distorting the shell elements, while the coupling of updated surface intersections is still enforced in the structural analysis stage. This relative movement of shell patches in the parametric coordinates helps maintain the quality of shell elements, even when intersections undergo significant relocation. Figure \ref{fig:shell-shopt-mint} shows the updated shape patches from the original design in Figure \ref{fig:shell-coupling}, where the parametric locations of the intersection are updated accordingly to reflect the relative movement of the shell patches during the optimization process. This method is achieved by considering the parametric locations of the surface intersection as state variables in the optimization framework and formulating a differentiable residual that accounts for the control points of shell patches with indices $k_1$ and $k_2$, along with the associated parametric locations of the intersection $\tilde{\bm{\xi}}^{k_1}$ and $\tilde{\bm{\xi}}^{k_2}$
\begin{align}
    \mathbf{R}_{\mathcal{L}} = \begin{bmatrix}
        \mathbf{N}^{k_1}(\tilde{\bm{\xi}}^{k_1})\mathbf{P}^{k_1} - \mathbf{N}^{k_2}(\tilde{\bm{\xi}}^{k_2})\mathbf{P}^{k_2} \\[4pt]
        h^{k_1}_{j} - h^{k_1}_{j-1} \\[4pt]
        \tilde{{\xi}}^{k_1}_a - 1\backslash 0\\[4pt]
        \tilde{{\xi}}^{k_1}_b - 1\backslash 0
    \end{bmatrix} = \mathbf{0} \text{ ,} \label{eq:intersection-residual}
\end{align}
where $h^{k_1}_{j}$ is the element size of the intersection in physical space. The first equation in \eqref{eq:intersection-residual} ensures that the locations of nodal points of the intersection, associated with two spline patches, coincide in physical space. The second equation in \eqref{eq:intersection-residual} ensures that adjacent physical elements at the intersection have the same length, thereby maintaining a uniform physical element size along the intersection. The last two equations specify the two end points of the intersection, which have parametric coordinates of 0 or 1 depending on the surface edge, assuming the spline patches have a unit square parametric domain and the lower left corner at (0,0).

The inclusion of parametric coordinates of patch intersections in the optimization problem impacts the structural analysis results, making the structural displacement dependent on both $\mathbf{P}$ and $\tilde{\bm{\xi}}$. This dependency is expressed as 
\begin{align}
    \mathbf{R}(\mathbf{P}, \tilde{\bm{\xi}}, \mathbf{d}) = \mathbf{0} \text{ .} \label{eq:nonmatching-residual-moving intersection}
\end{align}
The non-matching residual of structural analysis $\mathbf{R}$ and associated stiffness matrix $\mathbf{K} = \partial_{\mathbf{d}} \mathbf{R}$ are defined in \eqref{eq:nonmatching-residual-stiffness}. In contrast to the approach in Section \ref{subsubsec:shape-opt-ffd}, where the relative location between intersecting shell patches within one FFD-block remains fixed during the shape optimization process and $\tilde{\bm{\xi}}$ is not included in the optimization, this method requires the parametric locations of patch intersections to be treated as state variables. To enable gradient-based optimization, the partial derivatives $\partial_{\mathbf{P}} \mathbf{R}$ and $\partial_{\tilde{\bm{\xi}}} \mathbf{R}$ in \eqref{eq:nonmatching-residual-moving intersection} need to be computed. The derivation of these derivatives is discussed in detail in \cite[Section 3]{zhao2024shape}. It is noted that while the first-order derivatives of the spline basis functions are considered in the penalty energy for shell coupling, second-order derivatives are required in the partial derivative $\partial_{\tilde{\bm{\xi}}} \mathbf{R}$, thereby $C^1$ continuity across element boundaries is necessary. This requirement is naturally satisfied by the NURBS basis functions.

By incorporating differentiable intersections, the total derivative for the shape optimization problem is formulated as
\begin{align}
    \mathrm{d}_{\mathbf{P}} f = \partial_{\mathbf{P}} f - (\partial_{\mathbf{d}} f)^{\mathrm{T}} \mathbf{K}_{}^{-1} \left[\partial_{\mathbf{P}}\mathbf{R}_{} - \partial_{\tilde{\bm{\xi}}} \mathbf{R}_{} \, (\partial_{\tilde{\bm{\xi}}} \mathbf{R}_{\mathcal{L}})^{-1} \partial_{\mathbf{P}} \mathbf{R}_{\mathcal{L}}\right] \label{eq:derivative-multi-patch-df-dp} \text{ .}
\end{align}
Since the number of design variables is typically much larger than the number of model outputs, the adjoint method is employed to compute sensitivities in the optimization problem for enhanced efficiency.

\section{Design of GOLDFISH} \label{sec:design-goldfish}
Section \ref{subsec:software-dependencies} outlines the design of GOLDFISH and its software dependencies, which facilitate the open-source implementation. A discussion of the key OpenMDAO components for shape optimization of isogeometric shell structures is presented in Section \ref{subsec:optimization-components}.

\subsection{Software dependencies and workflow} \label{subsec:software-dependencies}

The design of GOLDFISH leverages the code generation capabilities in FEniCS to automate the computation of symbolic Gateaux derivatives for gradient-based optimization, while OpenMDAO is used to ensure modularity and flexibility across various design conditions and disciplines. The Python library is built on a suite of open-source software dependencies to streamline the entire design-analysis-optimization workflow.

The structural analysis is performed using PENGoLINS \cite{zhao2022open}, a Python library designed for complex shell structures modeled by isogeometric Kirchhoff--Love theory. In PNGoLINS, CAD geometries of shell structures are discretized isogeometrically and are directly available for analysis without FE mesh generation. Shell CAD geometries consisting of a collection of non-conforming NURBS surfaces are coupled using a penalty formulation \cite{herrema2019penalty}, where a penalty energy, as discussed in \eqref{eq:penalty-energy}, is integrated to preserve displacement continuity and angular compatibility along the intersection between shell patches. In the coupling procedure, a geometrically 2D, topologically 1D quadrature mesh is generated in the parametric space between two intersecting shell patches to serve as the integration domain for the penalty energy. This quadrature mesh is first positioned at the parametric location of the intersection with respect to the first shell patch to interpolate the displacement and covariant basis vectors. It then performs a similar operation with respect to the second shell patch. With the interpolated displacements and rotational quantities, the penalty energy and associated derivatives can be computed using FEniCS, enabling us to solve the coupled system of the complex shell structure. A schematic visualization of the coupling procedure is shown in Figure \ref{fig:shell-coupling}. Analysis results, including displacements and stress resultants for a range of benchmark problems and aerospace engineering applications, are presented in \cite{zhao2022open, zhao2024automated} to verify the accuracy of PENGoLINS.

PENGoLINS makes use of the Python interface of OpenCASCADE \cite{paviot2018pythonocc} as the geometry engine so that CAD geometries can be imported to extract knot vectors and control points directly and use them for IGA. A key feature inherited from OpenCASCADE is its ability to approximate NURBS patch intersections. This functionality is used to compute the parametric locations $\tilde{\bm{\xi}}$ of surface intersections, which define the integration domains for the penalty energy in \eqref{eq:penalty-energy}. Meanwhile, the computed parametric coordinates serve as the initial guess when solving the implicit equation between NURBS surface control points and intersection parametric coordinates in \eqref{eq:intersection-residual}. The IGA capabilities in PENGoLINS are powered by tIGAr, a Python library developed based on FEniCS that leverages the existing FE assembly routines. tIGAr constructs NURBS basis functions from Lagrange polynomials using the extraction technique \cite{BSEH11, Schillinger2016}, with the Galerkin approximation process fully automated using FEniCS.

\begin{figure}[!htb]\centering
    \includegraphics[width=0.9\textwidth]{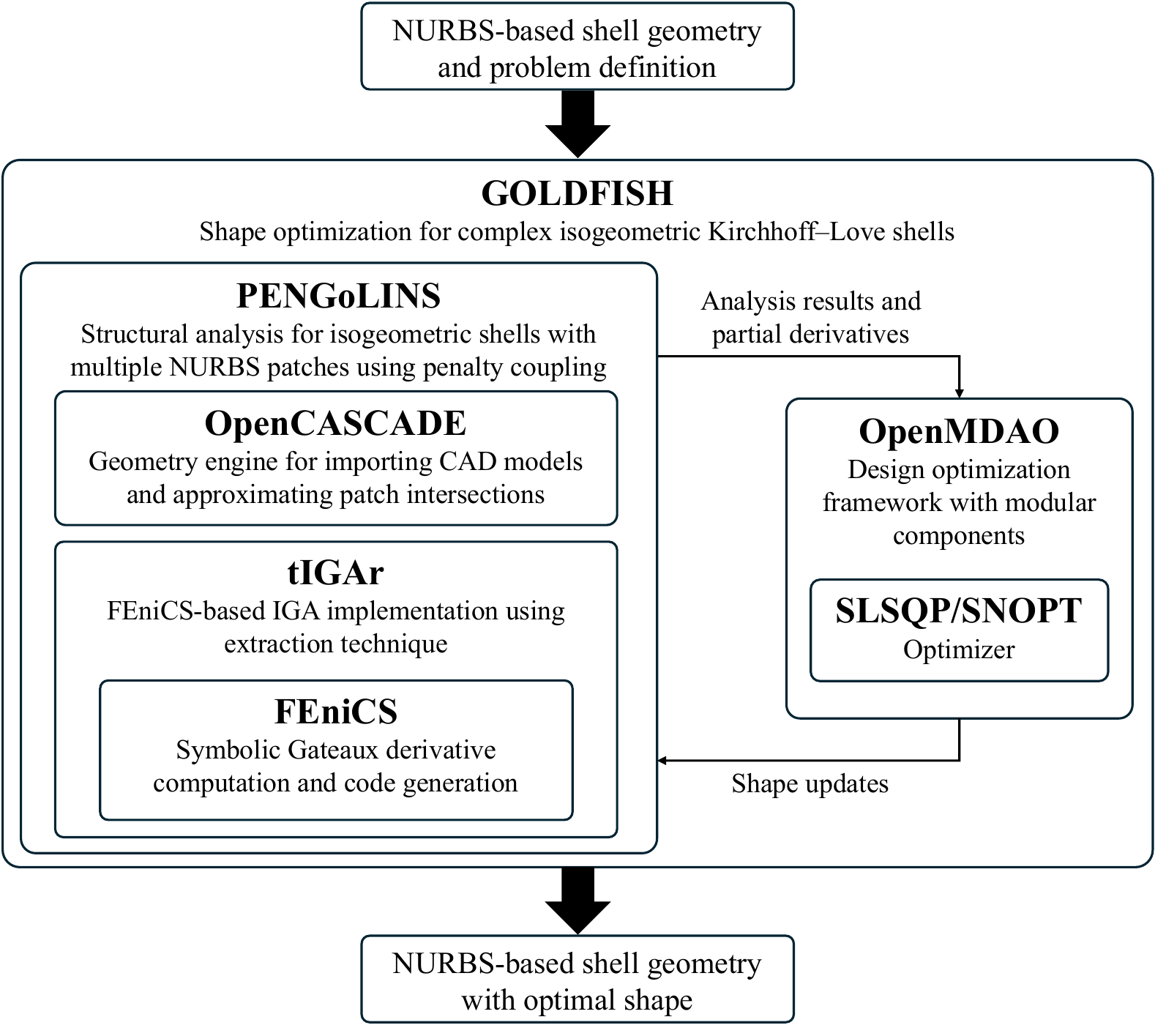}
    \caption{The design of the Python library GOLDFISH and its software dependencies. PENGoLINS is used for structural analysis and OpenMDAO is employed for numerical optimization. Analytical partial derivatives are computed in individual components in GOLDFISH. Both input and output of the software are NURBS-based shell geometries.}
    \label{fig:goldfish-code-structure}
\end{figure}

The numerical optimization is conducted using OpenMDAO \cite{gray2019openmdao}, which computes the total derivative of optimization problems using direct or adjoint methods, depending on the problem specifications. In shape optimization problems, the number of design variables is typically much larger than the number of model outputs, such as objective functions and constraints. OpenMDAO automatically organizes the partial derivatives provided by components into total derivatives using the adjoint method, significantly improving the efficiency of derivative computation. The modular design of OpenMDAO also facilitates and standardizes the implementation of individual components for the optimization problem. Each partial derivative in \eqref{eq:derivative-df-dpffd-full} and \eqref{eq:derivative-multi-patch-df-dp} can be implemented as a standard OpenMDAO component, which is connected automatically during the optimization process. This modular design greatly enhances the flexibility and applicability of the code framework, allowing it to be adapted to more customized problems. A series of essential components for shell shape optimization are discussed in Section \ref{subsec:optimization-components}. The optimization problem can be solved using the open-source optimizer SLSQP \cite{kraft1988software} or the commercial optimizer SNOPT \cite{gill2005snopt}. A schematic code structure of GOLDFISH is illustrated in Figure \ref{fig:goldfish-code-structure}.

Figure \ref{fig:goldfish-code-structure} illustrates the streamlined workflow of GOLDFISH, where users only need to provide the NURBS-based CAD geometry of the shell structure and define the optimization problem by specifying objective functions and constraints within OpenMDAO. The code framework then automatically performs the structural analysis and shape optimization on the NURBS-based geometry without FE mesh generation. A Lagrange extraction matrix is generated using \eqref{eq:lagrange-extraction-matrix-generation} for each shell patch to express the IGA spline bases by the Lagrange polynomial bases used in FEniCS. Both types of basis functions are associated with the 2D parametric meshes of the spline patches, as shown in Figure \ref{fig:shell-coupling}. These 2D parametric meshes are defined by the knot vectors and are generated automatically in FEniCS. The final output is also a NURBS-based geometry with updated control points that define the optimal shell shape. Throughout the optimization loop, shape updates and structural analysis are all performed directly on the CAD geometry. This integration considerably simplifies the shape optimization process.

\subsection{Optimization components of shell shape optimization} \label{subsec:optimization-components}
This section reviews the essential building blocks for an IGA-based shape optimization problem. For the shape optimization problem described in \eqref{eq:minmize-single-patch} with the associated total derivative \eqref{eq:derivative-df-dp-full}, the design variables are the coordinates of control points $\mathbf{P}$ of the NURBS surface defining the shell geometry. A standard InputsComp is created to provide the independent design variables to the following core components.
~\\
\begin{itemize}
    \item \textit{DispComp}: An implicit OpenMDAO component takes control points $\mathbf{P}$ of the shell surface as input and returns the corresponding displacement $\mathbf{d}$ by solving the isogeometric Kirchhoff--Love shell problem $\mathbf{R}_{\text{S}}(\mathbf{P}, \mathbf{d}) = \mathbf{0}$, where $\mathbf{R}_{\text{S}}$ is defined in \eqref{eq:KL-shell-residual}. Meanwhile, this component computes partial derivatives $\partial_{\mathbf{P}} \mathbf{R}_{\text{S}}$ and $\partial_{\mathbf{d}} \mathbf{R}_{\text{S}}$.\\
    
    \item \textit{ObjectiveComp}: An explicit OpenMDAO component calculates the objective function $f$ for the optimization problem, which typically depends on the shape of the shell and its displacements $f(\mathbf{P}, \mathbf{d})$. Additionally, this component provides partial derivatives $\partial_{\mathbf{P}} f$ and $\partial_{\mathbf{d}} f$.
\end{itemize}
~\\
With the partial derivatives computed by the two core components, the total derivative \eqref{eq:derivative-df-dp-full} is constructed automatically in OpenMDAO to guide shape updates until the optimal solution is reached.

\subsubsection{Components for FFD-based shape optimization} \label{subsubsec:components-ffd}
\begin{figure}[!htb]\centering
    \includegraphics[width=0.8\textwidth]{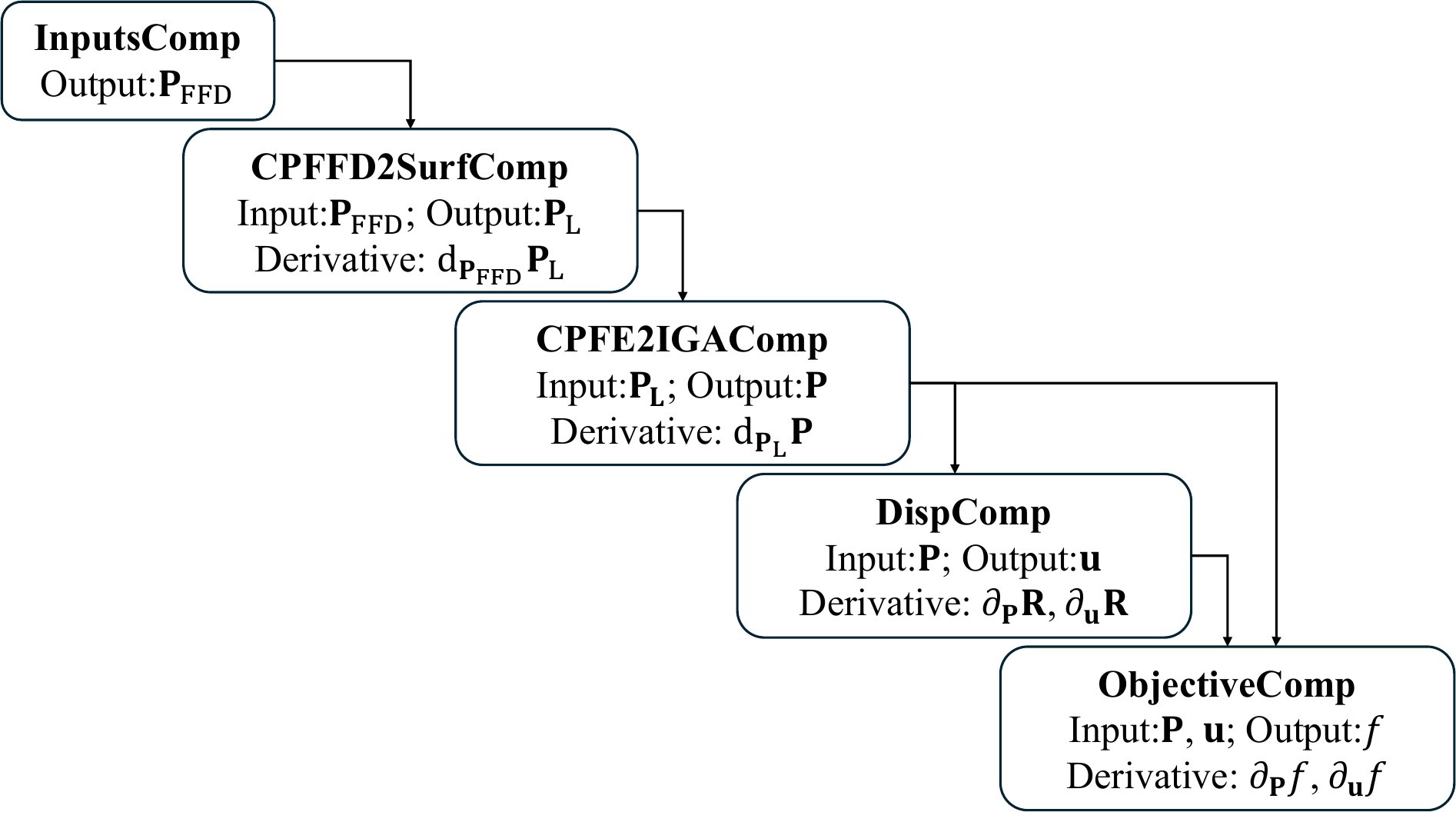}
    \caption{Component structure for shell shape optimization using the FFD-based approach.}
    \label{fig:shape-opt-structure-ffd}
\end{figure}

As discussed in Section \ref{subsubsec:shape-opt-ffd}, the FFD-based approach is applied to real-world CAD geometries consisting of multiple non-conforming NURBS patches to maintain the intersections. In this approach, control points of the trivariate B-spline block $\mathbf{P}_{\text{FFD}}$ serve as the design variables. Additional components are implemented for this approach to automate the optimization process, as outlined below.
~\\
\begin{itemize}
    \item \textit{CPFFD2SurfComp}: An explicit component computes the corresponding Lagrange nodal points $\mathbf{P}_{\text{L}}$ for the given control points of the FFD block $\mathbf{P}_{\text{FFD}}$, as illustrated in \eqref{eq:cp-ffd-cp-lagrange}--\eqref{eq:derivative-dpl-dpffd}. It also provides the derivative $\mathrm{d}_{\mathbf{P}_{\text{FFD}}} \mathbf{P}_{\text{L}}$ by evaluating the basis functions of the FFD block at the Lagrange nodal points of the shell surfaces in their initial configuration.\\

    \item \textit{CPFE2IGAComp}: An implicit component solves for the NURBS control points of the shell patches $\mathbf{P}$ given the input Lagrange nodal points $\mathbf{P}_{\text{L}}$, using the residual equation $\mathbf{M} \mathbf{P} - \mathbf{P}_{\text{L}} = \mathbf{0}$. The matrix $\mathbf{M}$ is the global extraction operator for the entire shell structure. This implicit component is employed to bypass the large matrix inversion as shown in \eqref{eq:derivative-dpk-dplk}. Since the degrees of freedom (DoFs) of $\mathbf{P}$ are fewer than $\mathbf{P}_{\text{L}}$, $\mathbf{M}$ is a nonsquare matrix. The resulting $\mathbf{P}$ is interpreted as a least-squares fit to $\mathbf{P}_{\text{L}}$. \\

    \item \textit{DispComp}: An implicit component solves for the displacement $\mathbf{d}$ of the multi-patch shell structure given the input shell NURBS control points $\mathbf{P}$. Unlike the shape optimization of single patch shell structure, the residual becomes $\mathbf{R}(\mathbf{P}, \mathbf{d}) = \mathbf{0}$ as discussed in \eqref{eq:nonmatching-residual-stiffness}, which includes a penalty energy to couple the intersecting shell patches. This component also returns the partial derivatives of the non-matching residual $\partial_{\mathbf{P}}\mathbf{R}$ and $\partial_{\mathbf{d}}\mathbf{R}$.
\end{itemize}
~\\
By connecting with the previously mentioned InputsComp and ObjectiveComp, we can perform shape optimization for the complex shell structures while preserving the non-matching patch intersection throughout the optimization process. Figure \ref{fig:shape-opt-structure-ffd} illustrates the component structure of the FFD-based shape optimization. This code structure is verified in the non-matching arch shape optimization example in Section \ref{subsec:example-arch}.

An illustrative implementation example is provided in the following code snippets, demonstrating the use of GOLDFISH within the Python environment. First, we import the OpenMDAO and GOLDFISH libraries.
\begin{lstlisting}[language=Python]
import openmdao.api as om
from GOLDFISH.nonmatching_opt_om import *
\end{lstlisting}
Next, a class \lstinline{ShapeOptGroupFFD} inherited from the OpenMDAO group is created for the FFD-based shape optimization problem, and the relevant parameters are initialized. The input of the class is an instance of the non-matching problem \lstinline{NonMatchingOptFFD}, which takes the CAD geometry, analysis definitions, and optimization conditions. These problem definitions are demonstrated in Section \ref{subsec:example-arch}.
\begin{lstlisting}[language=Python]
class ShapeOptGroupFFD(om.Group):
    def initialize(self):
        self.options.declare('nonmatching_opt_ffd')
    # Define optimization related parameters
    def init_parameters(self):
        self.nmopt_ffd = self.options['nonmatching_opt_ffd']
        self.opt_field = self.nmopt_ffd.opt_field
        self.init_cpffd_design = self.nmopt_ffd.\
                         shopt_init_cpffd_design
        self.input_cp_shapes = [cpffd.size for cpffd 
                                in self.init_cpffd_design]
\end{lstlisting}
A list of OpenMDAO components is then added to the group. The following code snippet shows the input component which takes the control points of the FFD block as the design variables. The \lstinline{CPFFD2SurfComp} and \lstinline{CPFE2IGAComp} components connect control points of the FFD block to the NURBS control points of shell patches using the formulation discussed in Section \ref{subsubsec:shape-opt-ffd}.
\begin{lstlisting}[language=Python]
    def setup(self):
        # Add inputs comp
        inputs_comp = om.IndepVarComp()
        for i, field in enumerate(self.opt_field):
            inputs_comp.add_output(
                VARNAME_CP_FFD_DESIGN+str(field),
                shape=self.input_cpffd_shapes[i],
                val=self.init_cpffd_design[i])
        self.add_subsystem('inputs_comp', inputs_comp, 
                           promotes=['*'])
        # Add FFD comp
        self.ffd2surf_comp = CPFFD2SurfComp(
            nonmatching_opt_ffd=self.nmopt_ffd)
        self.ffd2surf_comp.init_parameters()
        self.add_subsystem('CPFFD2Surf_comp', 
            self.ffd2surf_comp, promotes=['*'])
        # Add CPFE2IGA comp
        self.cpfe2iga_comp = CPFE2IGAComp(
            nonmatching_opt=self.nmopt_ffd)
        self.cpfe2iga_comp.init_parameters()
        self.add_subsystem('CPFE2IGA_comp', 
            self.cpfe2iga_comp, promotes=['*'])
\end{lstlisting}
Furthermore, the \lstinline{DispStatesComp} component performs the IGA on the updated CAD geometry and returns the structural response along with partial derivatives. In this example, we use the internal energy as the objective function, so the \lstinline{IntEnergyComp} component is added to the group.
\begin{lstlisting}[language=Python]
        # Add displacement comp
        self.disp_states_comp = DispStatesComp(
            nonmatching_opt=self.nmopt_ffd)
        self.disp_states_comp.init_parameters()
        self.add_subsystem('disp_comp', 
            self.disp_states_comp, promotes=['*'])
        # Add internal energy comp (obj function)
        self.int_energy_comp = IntEnergyComp(
            nonmatching_opt=self.nmopt_ffd)
        self.int_energy_comp.init_parameters()
        self.add_subsystem('int_energy_comp', 
            self.int_energy_comp, promotes=['*'])
\end{lstlisting}
Finally, we can specify the design variables and objective functions within the group to complete the setup of the optimization problem. Additionally, equality and inequality constraints can be specified in the group through the \lstinline{self.add_constraint} method.
\begin{lstlisting}[language=Python]
        # Add design variable and objective
        for i, field in enumerate(self.opt_field):
            self.add_design_var(
                VARNAME_CP_FFD_DESIGN+str(field),
                lower=DESVAR_L[i], upper=DESVAR_U[i])
        self.add_objective(VARNAME_INT_ENERGY)
\end{lstlisting}

\subsubsection{Components for moving intersections} \label{subsubsec:components-moving-intersection}
\begin{figure}[!htb]\centering
    \includegraphics[width=0.99\textwidth]{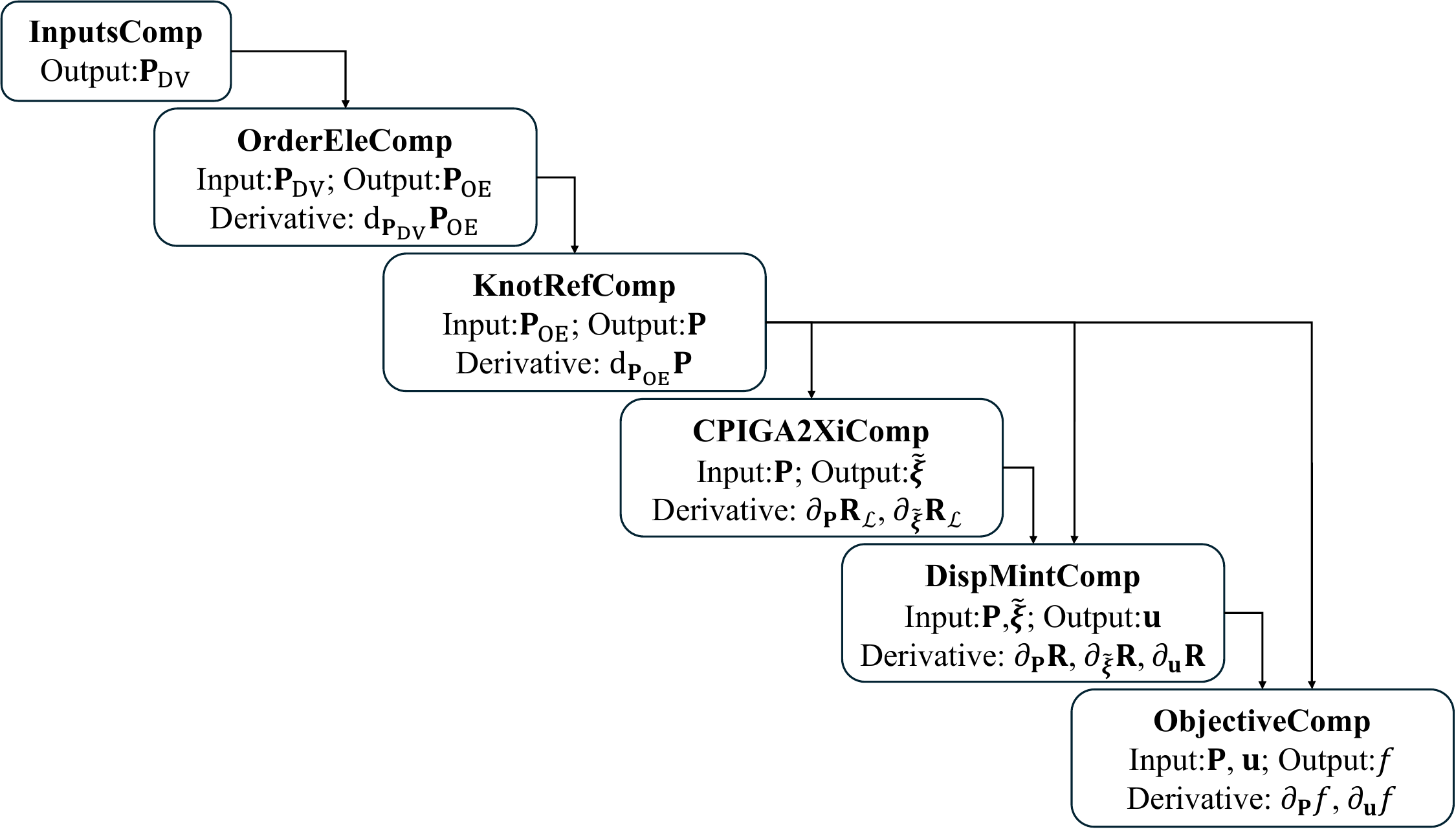}
    \caption{Component structure for shape optimization with moving intersections and the multilevel design method.}
    \label{fig:shape-opt-structure-multilevel}
\end{figure}

When the optimal shell structures require significant repositioning of intersections compared to the baseline design, we employ the moving intersections approach proposed in \cite{zhao2024shape}. This approach allows relative movement of shell patches, ensuring that shell elements maintain good quality in the optimized geometry. In this method, we implement the multilevel design approach \cite{nagy2010isogeometric, nagy2013isogeometric} to distinguish the design model and analysis model. As such, we can select the dimension of the design space independently from the analysis model, which typically has more DoFs for accurate analysis. By selecting a design model with much fewer DoFs than the analysis model, we can expedite the convergence of the optimizer while preserving the same geometry as the analysis model without introducing geometric errors. The multilevel design is achieved through order elevation and knot refinement of the NURBS surfaces in the design model. As a result, the design variables in this scenario are the control points of the design model. The essential components for shell shape optimization with moving intersections, combined with multilevel design, are listed below.
~\\
\begin{itemize}
    \item \textit{OrderElevationComp}: An explicit OpenMDAO component takes the control points of the coarse design model as input and returns corresponding control points after order elevation, following the relation $\mathbf{N}_{\text{DV}}(\bm{\xi}) \mathbf{P}_{\text{DV}} = \mathbf{N}_{\text{OE}}(\bm{\xi}) \mathbf{P}_{\text{OE}}$, where $\mathbf{N}_{\text{OE}}(\bm{\xi})$ represents higher-order NURBS basis functions than $\mathbf{N}_{\text{DV}}(\bm{\xi})$ but with the same unique knots. \\ %This relation can be expressed as a linear operation $ \mathbf{A}_{\text{OE}} \mathbf{P}_{\text{DV}} = \mathbf{P}_{\text{OE}}$, where the linear operator $\mathbf{A}_{\text{OE}}$ can be solved as $\mathbf{A}_{\text{OE}}= \left(\mathbf{N}_{\text{OE}}(\bm{\xi})^{\text{T}} \mathbf{N}_{\text{OE}}(\bm{\xi})\right)^{-1} \mathbf{N}_{\text{OE}}(\bm{\xi})^{\text{T}} \mathbf{N}_{\text{DV}}(\bm{\xi})$. Due to the number of DoFs in $\mathbf{P}_{\text{OE}}$, the matrix inversion can be easily computed to obtain $\mathbf{A}_{\text{OE}}$. \\
    
    \item \textit{KnotRefinementComp}: An explicit component computes control points of the analysis model $\mathbf{P}$ from the given input $\mathbf{P}_{\text{OE}}$ using knot refinement for NURBS surfaces. This is achieved using similar formulation $\mathbf{N}_{\text{OE}}(\bm{\xi}) \mathbf{P}_{\text{OE}} = \mathbf{N}_{}(\bm{\xi}) \mathbf{P}_{}$, where $\mathbf{N}_{}(\bm{\xi})$ has the same order as $\mathbf{N}_{\text{OE}}(\bm{\xi})$ but with refined interior knots. \\ %Due to the large size of $\mathbf{P}$, we construct the linear operator for knot refinement $\mathbf{A}_{\text{KR}}$ in an explicit way following the algorithm in \cite[Section 4]{piegl2012nurbs} to bypass the large matrix inversion. Thus, we can compute the updated control points of the analysis model through $\mathbf{A}_{\text{KR}} \mathbf{P}_{\text{OE}} = \mathbf{P}$ where $\mathbf{A}_{\text{KR}}$ is the associated derivative in this component. \\

    \item \textit{CPIGA2XiComp}: An implicit component takes control points of the analysis model $\mathbf{P}$ as inputs and solves the residual equation $\mathbf{R}_{\mathcal{L}}$ defined in \eqref{eq:intersection-residual} to determine the parametric coordinates of the movable intersections $\tilde{\bm{\xi}}$. Partial derivatives $\partial_{\mathbf{P}} \mathbf{R}_{\mathcal{L}}$ and $\partial_{\tilde{\bm{\xi}}} \mathbf{R}_{\mathcal{L}}$ are computed in this component.\\

    \item \textit{DispMintComp}: An implicit component similar to the \textit{DispComp} in Section \ref{subsubsec:components-moving-intersection} that solves structural displacement $\mathbf{d}$ using the penalty-based coupling formulation for isogeometric shells. In addition to control points of the shell patches $\mathbf{P}$, this component also takes the parametric coordinates of the moving intersections $\tilde{\bm{\xi}}$ as an input to formulate the residual $\mathbf{R}(\mathbf{P}, \tilde{\bm{\xi}}, \mathbf{d})$. This component also computes the partial derivatives $\partial_{\mathbf{P}}\mathbf{R}_{}$, $\partial_{\tilde{\bm{\xi}}}\mathbf{R}_{}$ and $\partial_{\mathbf{d}}\mathbf{R}_{}$.
\end{itemize}
~\\
The total derivative in \eqref{eq:derivative-multi-patch-df-dp}, along with a multilevel design method, are obtained by connecting the partial derivatives provided by the components mentioned above. The connections between individual components are outlined in Figure \ref{fig:shape-opt-structure-multilevel}. A numerical example verifies this code structure is demonstrated in Section \ref{subsec:example-tbeam}.

%Analogously, a T-beam geometry is demonstrated in Figure \ref{fig:shape-opt-structure-multilevel} to show the idea of multilevel design, where the black points are the control points of the coarse design model $\mathbf{P}_{\text{DV}}$ to modify the shape of the vertical patch, and the red and blue points are the control points in the analysis model $\mathbf{P}$ to improve the accuracy of the analysis. In this approach, we can opt to update the shape of the vertical patch only with relative movement between the two patches during the optimization process, and the element quality in the horizontal patch will not be affected. 

Code snippets are demonstrated to create the OpenMDAO group for the shape optimization of a complex shell structure with moving intersections. The class \lstinline{ShapeOptGroupMint} requires two arguments. The first is an instance of \lstinline{CPSurfDesign2Analysis} which provides the data and derivatives for the multilevel design of CAD geometry to reduce the dimension of the design space. The second argument is an instance of the non-matching problem \lstinline{NonMatchingOpt}, which is the parent class of \lstinline{NonMatchingOptFFD} without the FFD-related functions. Section \ref{subsec:example-tbeam} presents a numerical example using this group.
\begin{lstlisting}[language=Python]
class ShapeOptGroupMint(om.Group):
    def initialize(self):
        self.options.declare('cpdesign2analysis')
        self.options.declare('nonmatching_opt')
    def init_parameters(self):
        self.des2ana = self.options['cpdesign2analysis']
        self.nm_opt = self.options['nonmatching_opt']
        self.opt_field = self.nm_opt.opt_field
        self.init_cp_design = self.des2ana.init_cp_design
        self.input_cp_shapes = [len(cp) for cp 
                                in self.init_cp_design]
\end{lstlisting}
Next, we add the input component which takes control points of the coarse CAD geometry as the design variables. The order elevation and the knot refinement components are included to perform the $k$-refinement for the multilevel design, producing the fine analysis model.
\begin{lstlisting}[language=Python]
    def setup(self):
        # Add inputs comp
        inputs_comp = om.IndepVarComp()
        for i, field in enumerate(self.opt_field):
            inputs_comp.add_output(
                VARNAME_CP_SURF_COARSE+str(field),
                shape=self.input_cp_shapes[i],
                val=self.init_cp_design[i])
        self.add_subsystem('input_comp', inputs_comp, 
                           promotes=['*'])
        # Add order elevation comp
        self.cp_order_ele_comp = CPSurfOrderElevationComp(
            cpdesign2analysis=self.des2ana)
        self.cp_order_ele_comp.init_parameters()
        self.add_subsystem('CP_order_ele_comp', 
            self.cp_order_ele_comp, promotes=['*'])
        # Add knot refinement comp
        self.cp_knot_refine_comp = CPSurfKnotRefinementComp(
            cpdesign2analysis=self.des2ana)
        self.cp_knot_refine_comp.init_parameters()
        self.add_subsystem('CP_knot_refine_comp', 
            self.cp_knot_refine_comp, promotes=['*'])
\end{lstlisting}
Subsequently, the \lstinline{CPIGA2XiComp} component is added to calculate the parametric coordinates of the intersections for a given set of surface control points by solving the implicit equation \eqref{eq:intersection-residual} and computing the partial derivatives. The displacement component for moving intersections \lstinline{DispMintComp} is then added to evaluate the structural responses for updated surface control points and intersection locations. Similarly, the internal energy component is used to define the objective function.
\begin{lstlisting}[language=Python]
        # Add CPIGA2Xi comp
        self.cpiga2xi_comp = CPIGA2XiComp(
            nonmatching_opt=self.nm_opt)
        self.cpiga2xi_comp.init_parameters()
        self.add_subsystem('CPIGA2xi_comp', 
            self.cpiga2xi_comp, promotes=['*'])
        # Add displacement comp with moving int
        self.disp_states_comp = DispMintStatesComp(
            nonmatching_opt=self.nm_opt)
        self.disp_states_comp.init_parameters()
        self.add_subsystem('disp_comp', 
            self.disp_states_comp, promotes=['*'])
        # Add internal energy comp (objective function)
        self.int_energy_comp = IntEnergyComp(
            nonmatching_opt=self.nm_opt)
        self.int_energy_comp.init_parameters()
        self.add_subsystem('int_energy_comp', 
            self.int_energy_comp, promotes=['*'])
\end{lstlisting}
Lastly, we assign the coarse control points of the shell patches as design variables and designate the internal energy as the objective function to complete the problem setup.
\begin{lstlisting}[language=Python]
        for i, field in enumerate(self.opt_field):
            self.add_design_var(
                VARNAME_CP_SURF_COARSE+str(field),
                lower=DESVAR_L[i], upper=DESVAR_U[i])
        self.add_objective(VARNAME_INT_ENERGY)
\end{lstlisting}
The preprocessing and setup of the non-matching problems, as well as the usage of the OpenMDAO groups mentioned above, are discussed in detail in Section \ref{sec:numerical-examples}.

\section{Numerical examples} \label{sec:numerical-examples}
In this section, we present the implementation of benchmark problems to demonstrate the use of GOLDFISH and validate it with reference solutions. Then we showcase the application of GOLDFISH to the design optimization of aircraft wings.

\subsection{Non-matching arch shape optimization}\label{subsec:example-arch}
\begin{figure}[!htbp]
    \centering
    \begin{subfigure}[!htbp]{0.49\textwidth}
        \includegraphics[width=\textwidth]{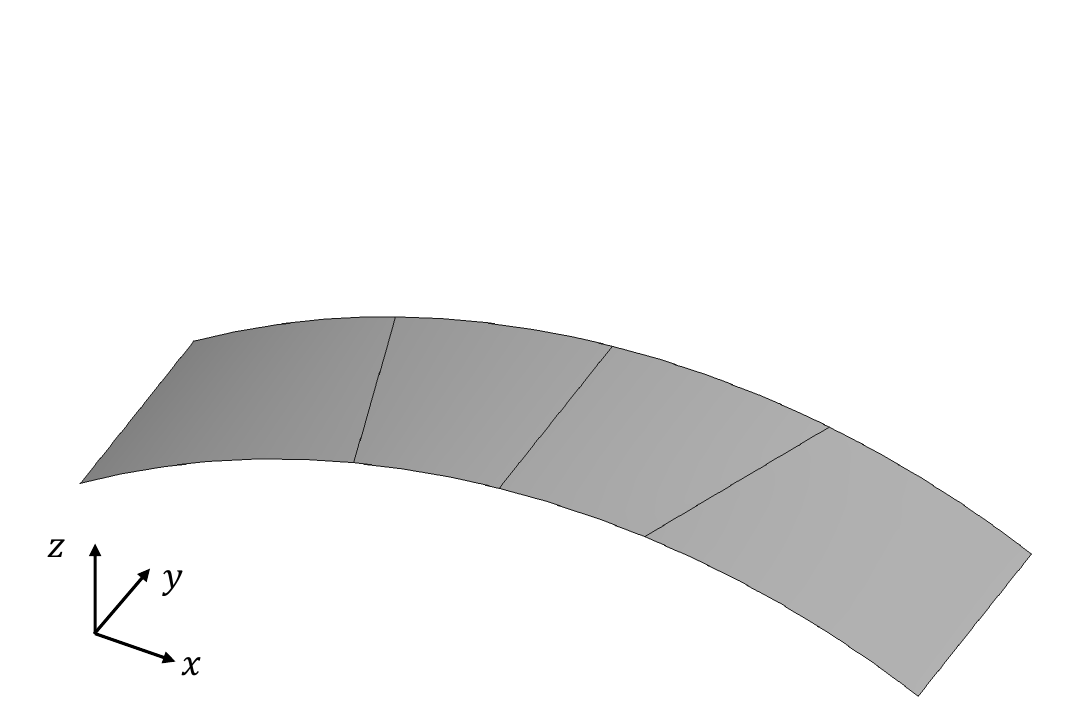}
        \caption{}
        \label{subfig:arch-geom-init}
    \end{subfigure}
    \hfill
    \begin{subfigure}[!htbp]{0.49\textwidth}
        \includegraphics[width=\textwidth]{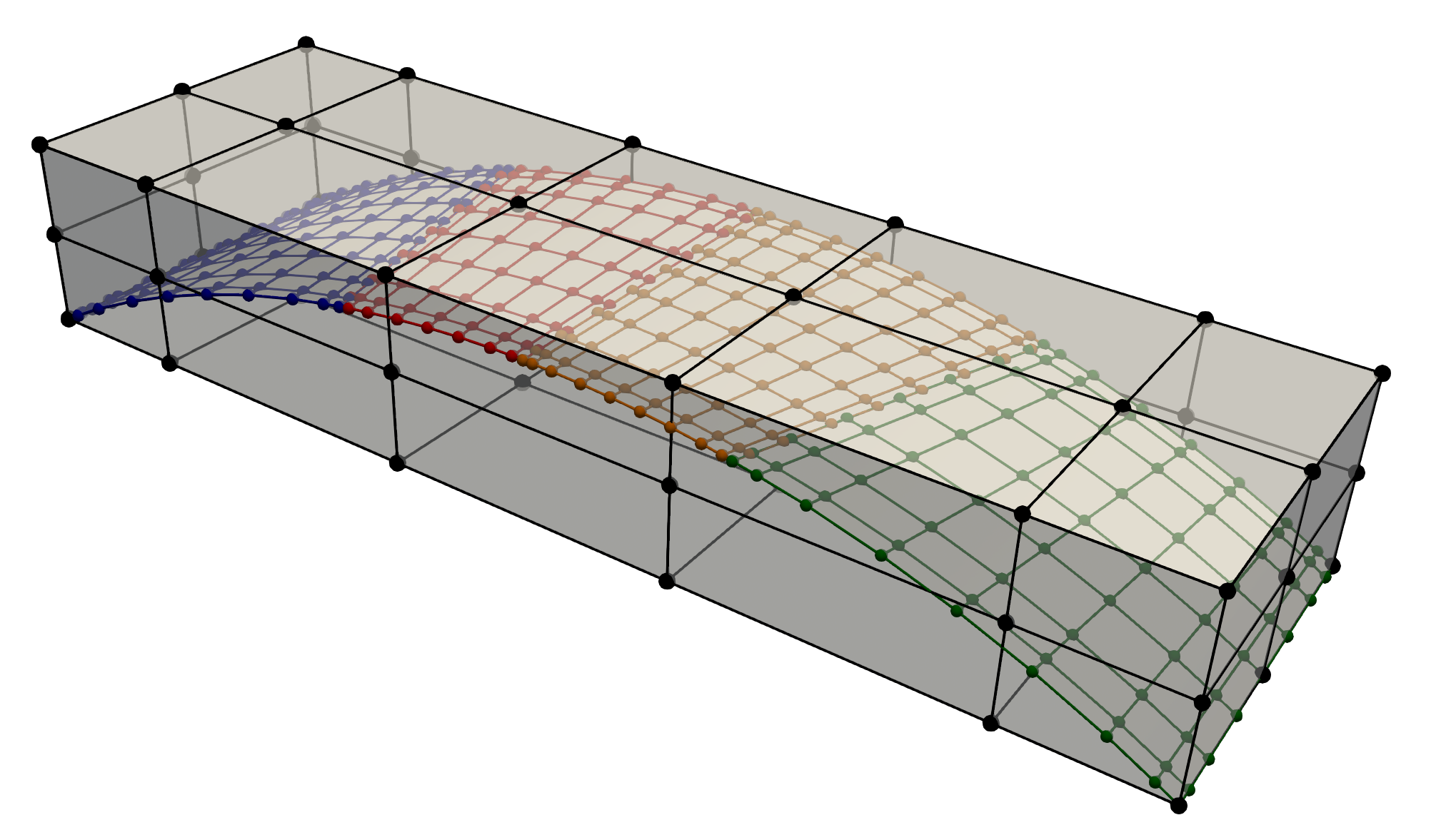}
        \caption{}
        \label{subfig:arch-geom-ffd}
    \end{subfigure}
    \caption{(a) The initial design of an arch geometry consisting of four non-matching NURBS patches. (b) The initial arch geometry is embedded in a 3D B-spline block for FFD-based shape optimization.}
    \label{fig:arch-geom}
\end{figure}

We use the arch shape optimization problem to verify the accuracy of the FFD-based shape update scheme in GOLDFISH. The benchmark problem was first proposed in \cite[Section 8]{kiendl2014isogeometric}, where an arch geometry is subjected to a distributed downward load and fixed at both edges. The optimal shape of the arch with minimum internal energy is a quadratic parabola with an analytical height-to-length ratio so that the external load is entirely supported by membrane forces \cite[Section 8]{kiendl2014isogeometric}. To examine the capability of the FFD-based approach to preserve non-conforming surface intersections during shape updates, we create an arch geometry consisting of four non-conforming NURBS patches, as shown in Figure \ref{subfig:arch-geom-init}. The arch geometry has a width of 3 m, a length of 10 m, and a shell thickness of 0.01 m. Young's modulus of 1000 GPa and Poisson's ratio of 0.0 are used for material properties. A 3D B-spline FFD block is generated to enclose the initial arch geometry, with the isogeometric discretization of both the FFD block and the arch geometry demonstrated in Figure \ref{subfig:arch-geom-ffd}.
The following listings provide essential implementation details for the FFD-based shape optimization using GOLDFISH.

% Listing \ref{lst:Python-listing}
% \begin{lstlisting}[language=Python, caption={Listing test}, label={lst:Python-listing}]
We first define the material properties, coupling coefficient, and other basic parameters. Considering the three physical directions correspond to $[0,1,2]$ in implementation since Python is a zero-based indexing language, we select control points in the $z$ direction for optimization and define \lstinline{opt_field} as $[2]$.

\begin{lstlisting}[language=Python]
E = Constant(1.0e12) # Young's modulus, Pa
nu = Constant(0.) # Poisson's ratio
h_th = Constant(0.01) # Shell thickness, m
pressure = Constant(1.) # Pressure magnitude, Pa
penalty_coefficient = 1.0e3 # Penalty coefficient
opt_field = [2] # Optimize z-coordinates only
ffd_block_nel = [4,1,1] # Nel of FFD block
p_ffd = 2 # Degree of the B-spline block
\end{lstlisting}
To perform the FFD-based shape optimization for the non-matching shells, we import the initial CAD geometry of the arch in IGES or STEP format into the running process using the Python interface of OpenCASCADE.
\begin{lstlisting}[language=Python]
topo_shapes = read_igs_file("init_arch_geom.igs", 
                           as_compound=False)
# Surface type conversion
occ_surf_list = [topoface2surface(face, BSpline=True) 
                 for face in topo_shapes] 
num_surfs = len(occ_surf_list)
\end{lstlisting}
We then create a geometry preprocessor instance to find all patch intersections and compute their parametric coordinates.
\begin{lstlisting}[language=Python]
preproc = OCCPreprocessing(occ_surf_list)
preproc.compute_intersections(mortar_refine=2)
\end{lstlisting}
Next, we generate a list of tIGAr spline instances to build the extraction matrices. These matrices are utilized in the FFD-based shape update and IGA. The implementation of function \lstinline{OCCBSpline2tIGArSpline}, which is standard to create a tIGAr spline instance from a given spline surface containing the knot vectors and control points, can be found in the GOLDFISH repository. Fixed boundary conditions are applied to the first and last shell patches and are implemented in this function. We omit these details here to focus on the setup of FFD-based shape optimization.
\begin{lstlisting}[language=Python]
splines = []
for i in range(num_surfs):
    splines += [OCCBSpline2tIGArSpline(
                preproc.BSpline_surfs[i])]
\end{lstlisting}
The next step is to create the non-matching coupling instance for the list of tIGAr spline instances using \lstinline{NonMatchingOptFFD}. This allows us to perform automated IGA and set up the FFD-based shape optimization problem.
\begin{lstlisting}[language=Python]
nmopt_ffd = NonMatchingOptFFD(splines, E, h_th, nu)
nmopt_ffd.create_mortar_meshes(preproc.mortar_nels)
# Optimize z-coords for all shell patches
nmopt_ffd.set_shopt_surf_inds_FFD(opt_field, [0,1,2,3])
# Quadrature meshes setup for penalty energy
nmopt_ffd.mortar_meshes_setup(preproc.mapping_list, 
          preproc.intersections_para_coords, 
          penalty_coefficient)
\end{lstlisting}
A 3D B-spline FFD block is created by specifying the number of elements, degrees, and limits of the control points in the three directions using function \lstinline{create_3D_block}. 
\begin{lstlisting}[language=Python]
# Create the 3D B-spline FFD block
cp_lims = nonmatching_opt_ffd.cpsurf_des_lims
for field in opt_field:
    cp_range = cp_lims[field][1]-cp_lims[field][0]
    cp_lims[field][1] = cp_lims[field][1]+0.2*cp_range
FFD_block = create_3D_block(ffd_block_nel, p_ffd, cp_lims)
\end{lstlisting}
By providing the knot vectors and control points of the FFD block to the non-matching coupling instance and arranging the related constraints on the FFD block control points, we complete the setup of the FFD block for shape optimization. The method \lstinline{set_shopt_align_CPFFD} eliminates redundant design variables by aligning control points in the width direction. The \lstinline{set_shopt_pin_CPFFD} method fixes control points on the lower edges of the FFD block, while \lstinline{set_shopt_regu_CPFFD} ensures control points stay within the range of their adjacent neighbors, preventing unrealistic shapes. 
\begin{lstlisting}[language=Python]
# Set FFD block to the shell optimization problem
nmopt_ffd.set_shopt_FFD(FFD_block.knots, FFD_block.control)
# Set CP alignment in the width direction
nmopt_ffd.set_shopt_align_CPFFD(align_dir=[[1]])
# Set constraints to fix two lower edges of the block
nmopt_ffd.set_shopt_pin_CPFFD(pin_dir0=[2], pin_side0=[[0]], 
                              pin_dir1=[1], pin_side1=[[0]])
# Set constraints to prevent self-penetration
nmopt_ffd.set_shopt_regu_CPFFD()
\end{lstlisting}
A list of the PDE residual forms based on the Kirchhoff--Love shell theory with St. Venant Kirchhoff material model is generated for all shell patches. These residual forms are expressed as FEniCS Unified Form Language (UFL) \cite{Alnaes2014} Form objects, which allows for automatic computation of symbolic derivatives. By inputting the residual forms into \lstinline{nmopt_ffd}, the stiffness matrix $\mathbf{K}$ and residual vector $\mathbf{R}$ of the non-matching shell structure are assembled in PENGoLINS subroutines.
\begin{lstlisting}[language=Python]
source_terms = []
residuals = []
for i in range(num_surfs):
    X = nmopt_ffd.splines[i].F # Shell geometry
    # Calculate curvilinear basis vectors
    A0,A1,A2,_,_,_ = surfaceGeometry(nmopt_ffd.splines[i], X)
    v_vec = as_vector([Constant(0.), Constant(0.), 
                       Constant(1.)])
    # Constant downward distributed pressure
    force = as_vector([Constant(0.), Constant(0.), 
                       -pressure*inner(v_vec, A2)])
    source_terms += [inner(force, nmopt_ffd.splines[i]\
        .rationalize(nmopt_ffd.spline_test_funcs[i]))\
        *nmopt_ffd.splines[i].dx]
    residuals += [SVK_residual(nmopt_ffd.splines[i], 
                  nmopt_ffd.spline_funcs[i], 
                  nmopt_ffd.spline_test_funcs[i], 
                  E, nu, h_th, source_terms[i])]
nmopt_ffd.set_residuals(residuals)
\end{lstlisting}
Subsequently, we construct the FFD-based shape optimization model using the component structure and implementation discussed in Section \ref{subsubsec:shape-opt-ffd}, and then create the OpenMDAO optimization problem.
\begin{lstlisting}[language=Python]
# Create the FFD-based shape optimization model
model = ShapeOptGroupFFD(nonmatching_opt_ffd=nmopt_ffd)
model.init_parameters()
prob = om.Problem(model=model)
\end{lstlisting}
Finally, the SLSQP optimizer with a tolerance of $10^{-12}$ is selected for this problem, and the objective function is minimized using method \lstinline{run_driver()}.
\begin{lstlisting}[language=Python]
prob.driver = om.ScipyOptimizeDriver()
prob.driver.options['optimizer'] = 'SLSQP'
prob.driver.options['tol'] = 1e-12
prob.driver.options['maxiter'] = 1000
# Set up and run the optimization problem
prob.setup()
prob.run_driver()
\end{lstlisting}

\begin{figure}[!htbp]
    \centering
    \begin{subfigure}[!htbp]{0.49\textwidth}
        \includegraphics[width=\textwidth]{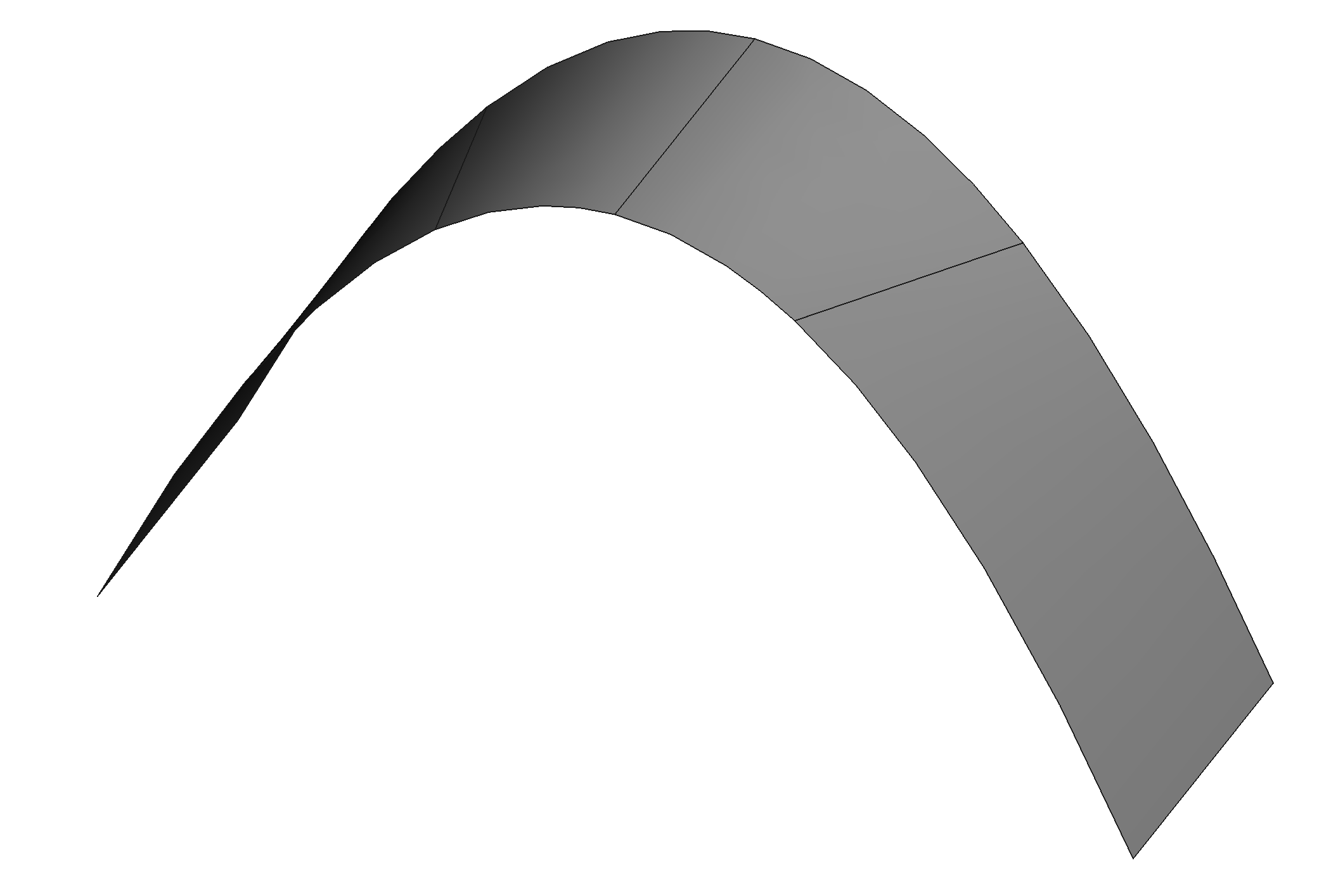}
         \caption{}
        \label{subfig:arch-geom-opt}
    \end{subfigure}
    \hfill
    \begin{subfigure}[!htbp]{0.49\textwidth}
        \includegraphics[width=\textwidth]{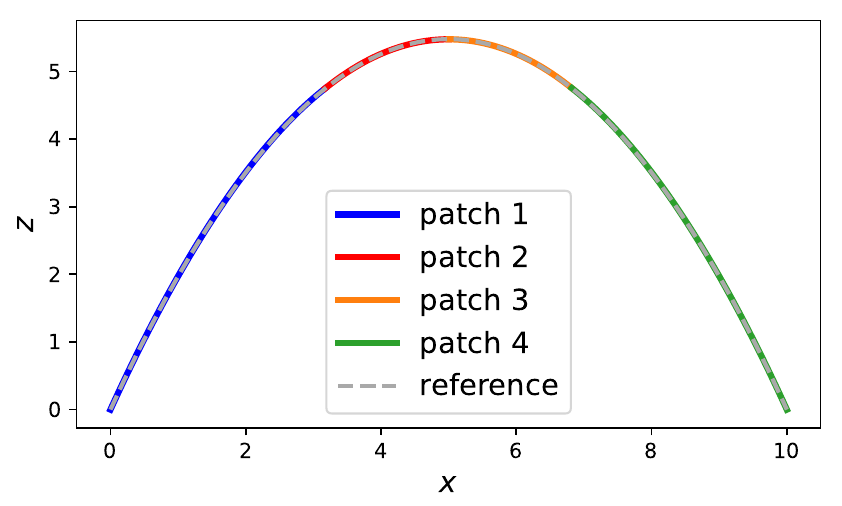}
        \caption{}
        \label{subfig:arch-opt-slice}
    \end{subfigure}
    \caption{(a) The optimized design of the arch geometry. (b) The cross-sectional view of the optimized arch compared with the analytical optimum.}
    \label{fig:arch-opt}
\end{figure}

The optimized arch geometry after 40 iterations is shown in Figure \ref{subfig:arch-geom-opt}. A comparison of the sliced view of the optimized arch with the analytical optimal shape is presented in Figure \ref{subfig:arch-opt-slice}, where the optimized geometry closely aligns with the reference solution from \cite{kiendl2014isogeometric}. The height-to-length ratio of the optimized solution is 5.4748, which shows a negligible difference from the analytical value of 5.4779.

\subsection{T-beam shape optimization with moving intersections} \label{subsec:example-tbeam}
\begin{figure}[!htbp]
    \centering
    \begin{subfigure}[!htbp]{0.49\textwidth}
        \includegraphics[width=\textwidth]{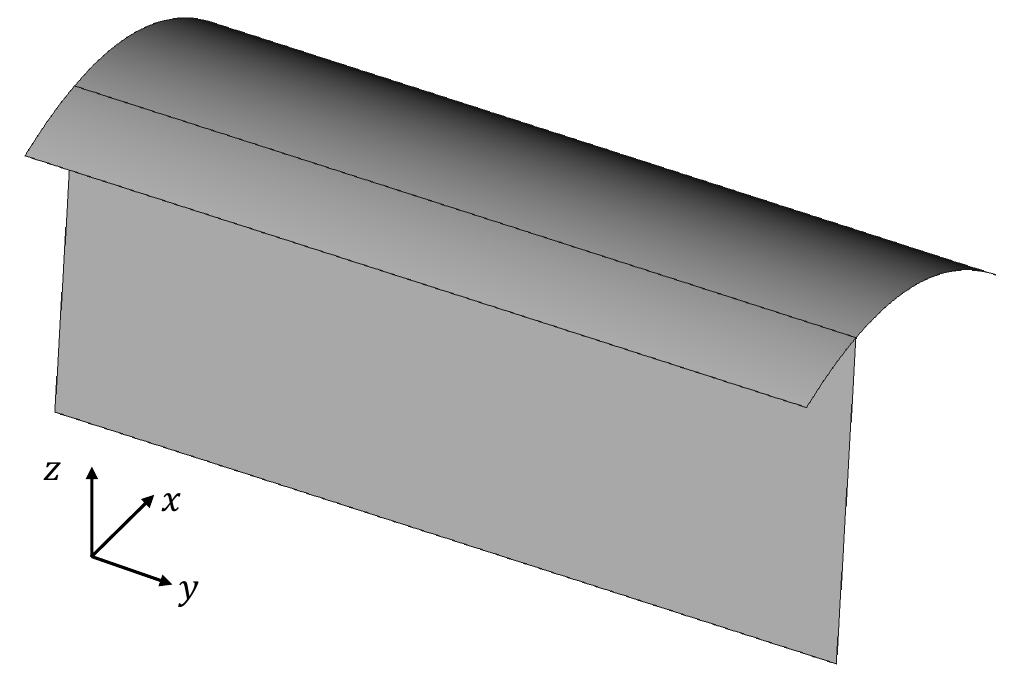}
        \caption{}
        \label{subfig:tbeam-geom-init}
    \end{subfigure}
    \hfill
    \begin{subfigure}[!htbp]{0.49\textwidth}
        \includegraphics[width=\textwidth]{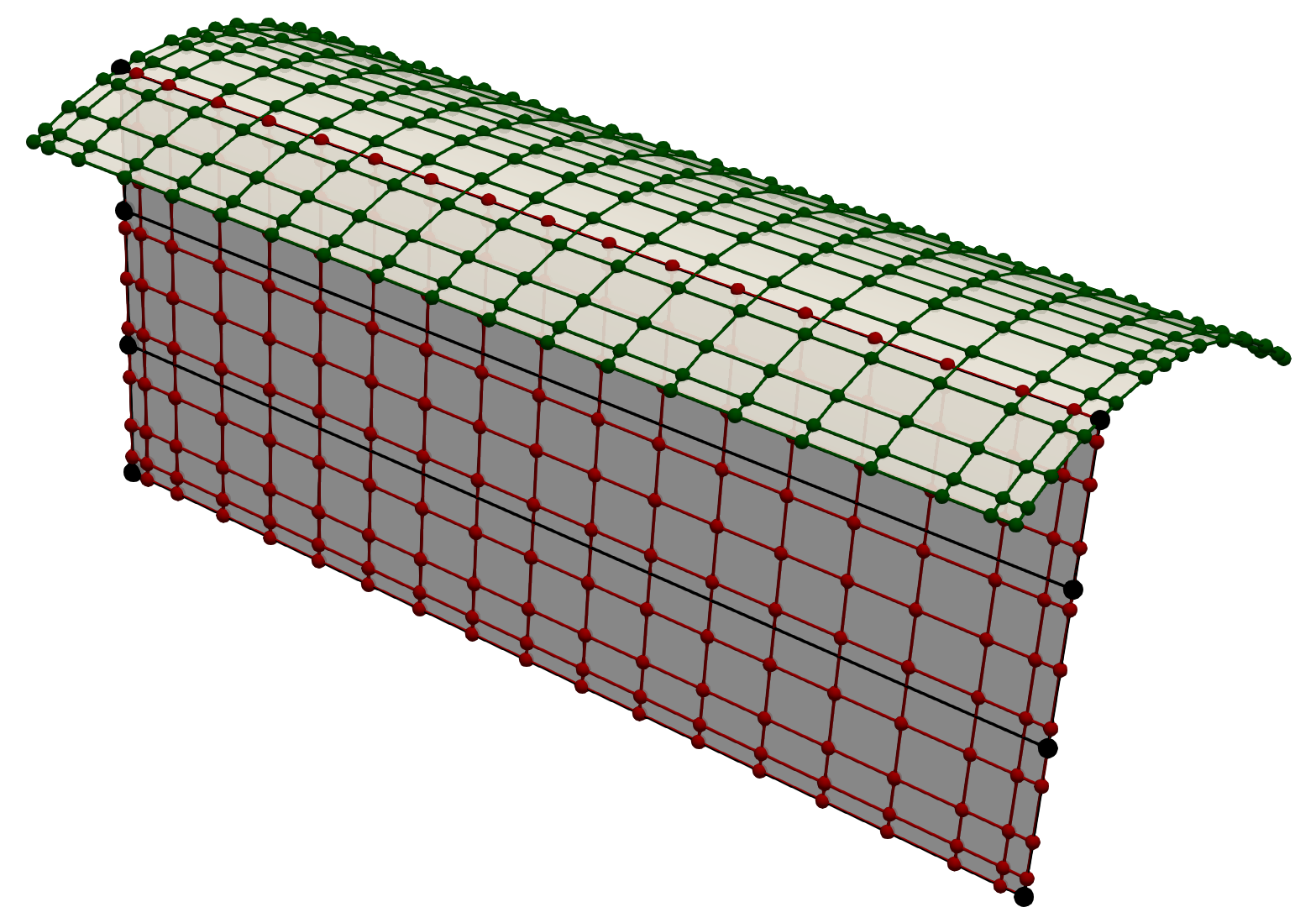}
        \caption{}
        \label{subfig:tbeam-geom-multilevel}
    \end{subfigure}
    \caption{(a) The baseline design of the T-beam geometry consists of two B-spline patches. The top patch is a curved surface. (b) Isogeometric discretization of the initial T-beam, where the red and green control points represent DoFs in the analysis model, and black control points represent DoFs in the design model for shape optimization.}
    \label{fig:tbeam-init}
\end{figure}

In this section, we use a T-beam geometry to demonstrate the use of GOLDFISH for shape optimization of multi-patch shell structures with moving intersections. The initial T-beam geometry is described by two B-spline surfaces. The top patch is a parabolic curved surface with dimensions of 2 m in width, 5 m in length, and spans from 0 m to 0.3 m in height. The vertical patch is a flat surface with dimensions of 2 m $\times$ 5 m and is positioned at one quarter of the top patch in the horizontal direction in the baseline design. The CAD geometry of the T-beam, shown in Figure \ref{subfig:tbeam-geom-init}, is subjected to a distributed pressure in the downward vertical direction and is fixed at one end. 

The optimal design with minimum internal energy is obtained when the vertical patch is positioned at the center of the top patch, maintaining a constant volume. In this example, we optimize the shape of the vertical patch while keeping the geometry of the top patch fixed. Meanwhile, the intersection between the two patches is allowed to move during the optimization process. Both shell patches have a thickness of 0.1 m, with material properties of Young's modulus $E=10^7$ Pa and Possion's ratio $\nu=0$. The discretization of the coarse design model for the vertical patch is indicated by black lines in Figure \ref{subfig:tbeam-geom-multilevel}, while the red and green points denote the discretization of the fine analysis model. The problem setup and GOLDFISH implementation are illustrated in the following code snippets.

\begin{figure}[!htbp]
    \centering
    \begin{subfigure}[!htbp]{0.49\textwidth}
        \includegraphics[width=\textwidth]{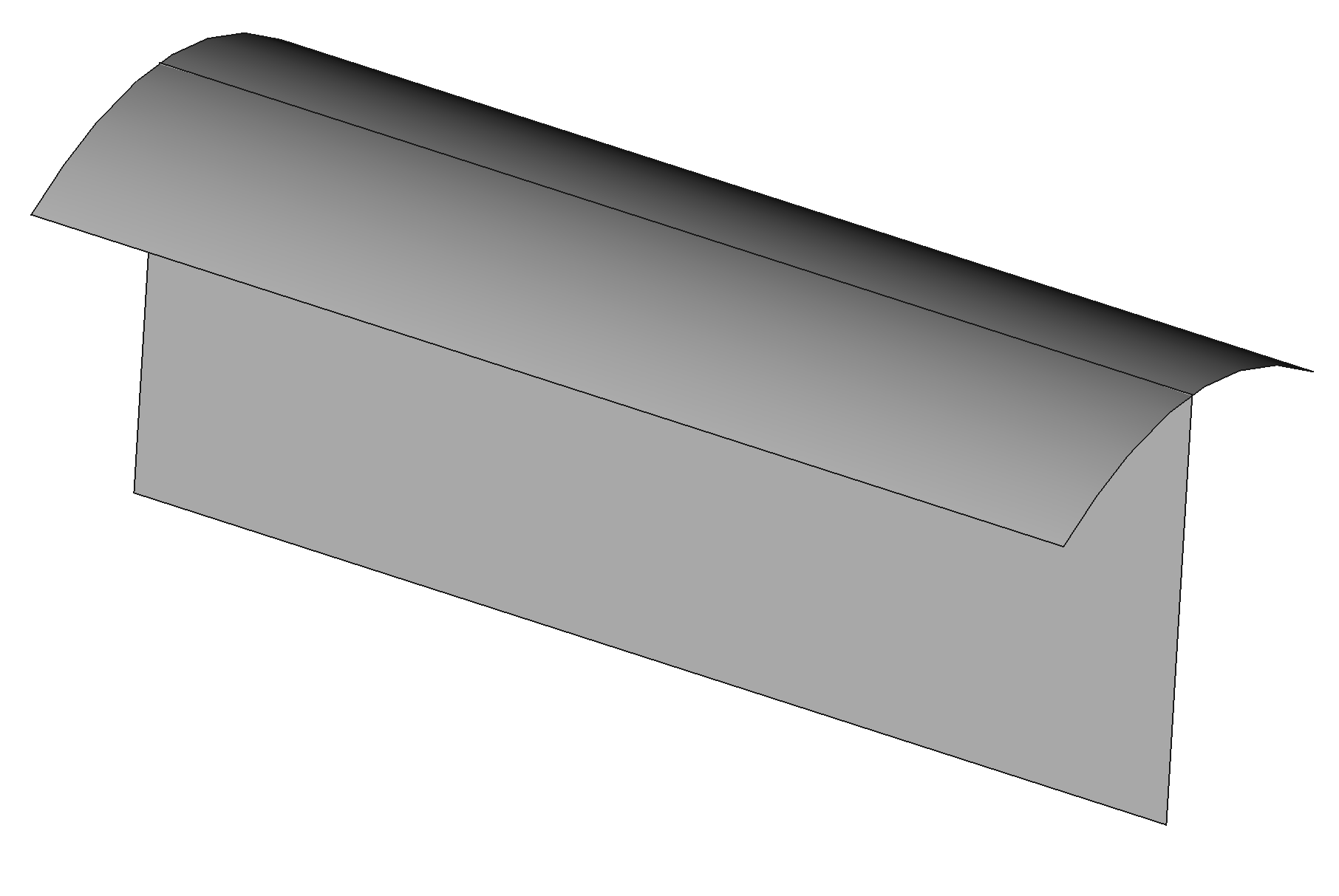}
        \caption{}
        \label{subfig:tbeam-geom-opt}
    \end{subfigure}
    \hfill
    \begin{subfigure}[!htbp]{0.49\textwidth}
        \includegraphics[width=\textwidth]{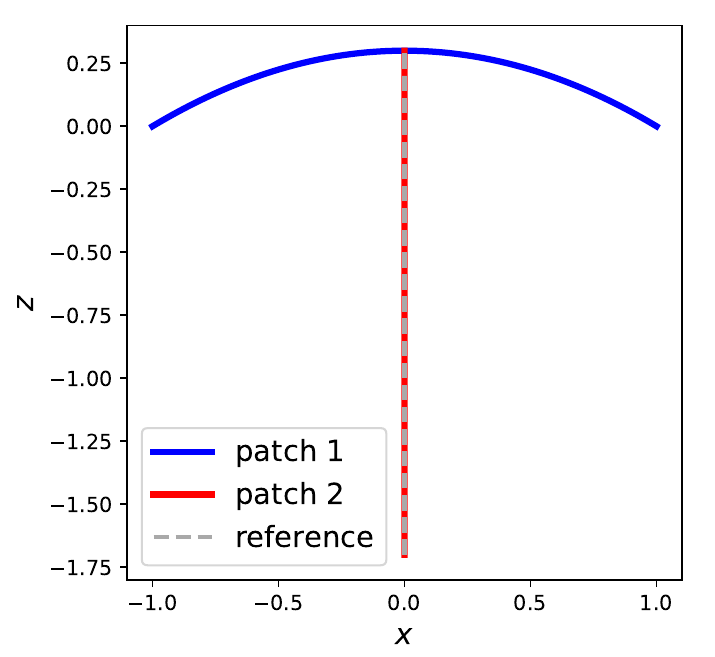}
        \caption{}
        \label{subfig:tbeam-opt-slice}
    \end{subfigure}
    \caption{(a) The optimized T-beam geometry with a curved top patch. (b) Cross-sectional view of the optimized T-beam, the vertical patch is moved to the center of the top patch and maintains the T-junction.}
    \label{fig:tbeam-opt}
\end{figure}

The problem parameters definitions and CAD geometry import are similar to the previous example and will not be repeated. For the shape optimization, an instance of the geometry processor \lstinline{OCCPreprocessing} and an instance of the non-matching problem \lstinline{NonMatchingOpt} are created.
% \begin{lstlisting}[language=Python]
% preproc = OCCPreprocessing(occ_surf_list)
% preproc.compute_intersections(mortar_refine=2)
% nmopt = NonMatchingOpt(splines, E, h_th, nu)
% \end{lstlisting}
We specify the fields of the control points for optimization as $x$ and $z$ coordinates, corresponding to \lstinline{opt_field} as $[0,2]$. Further, we specify the surface indices to be optimized in each field. For the vertical patch with an index of 1, the shape optimization surface indices are $[[1], [1]]$. We proceed to create and set up the quadrature meshes using \lstinline{mortar_meshes_setup}, similar to the previous example. The key difference is that the argument \lstinline{transfer_mat_deriv}, which defaults to 1, is set to 2 in this case since the partial derivative of the non-matching residual with respect to intersection parametric coordinates requires second-order derivatives of the spline basis functions. %, as discussed in Section \ref{subsubsec:shape-opt-moving-int}.
\begin{lstlisting}[language=Python]
opt_field = [0,2]
shopt_surf_inds = [[1], [1]]
nmopt.set_shopt_surf_inds(opt_field, shopt_surf_inds)
nmopt.set_geom_preprocessor(preproc)
nmopt.create_mortar_meshes(preproc.mortar_nels)
nmopt.mortar_meshes_setup(preproc.mapping_list, 
      preproc.intersections_para_coords, 
      penalty_coefficient, transfer_mat_deriv=2)
\end{lstlisting}
We use the \lstinline{check_intersections_type} method to check the types of intersections. An intersection is treated as differentiable if it is not located at the edges of both intersecting spline patches. Next, the method \lstinline{create_diff_intersections} is used to generate associated data for the differentiable intersections. The argument \lstinline{num_edge_pts} is a list of integers, where each item specifies the number of points in the quadrature mesh used to enforce the T-junction. For this T-beam example, there is only one T-junction, and since the vertical patch is set to remain straight in the axial direction during optimization, a single point is sufficient to ensure the T-junction.
\begin{lstlisting}[language=Python]
preproc.check_intersections_type()
preproc.get_diff_intersections()
nmopt.create_diff_intersections(num_edge_pts=[1])
\end{lstlisting}
We then use the geometry preprocessor to initialize an instance of \lstinline{CPSurfDesign2Analysis} to establish the multilevel design framework between the design model and the analysis model. 
\begin{lstlisting}[language=Python]
des2ana = CPSurfDesign2Analysis(preproc, opt_field, 
                                shopt_surf_inds)
\end{lstlisting}
We assume the imported CAD geometry represents the analysis model in the workflow. Therefore, we need to define the space for the design model. We first specify the order and knots of the design model. In this example, the horizontal location of the vertical patch is described by a cubic B-spline in the $\xi_1$ direction and a linear B-spline in the $\xi_2$ direction, while the vertical location is described by linear B-splines in both directions. All spline curves in the design model have a single knot span. The orders and knot vectors for the design model are defined as follows.
\begin{lstlisting}[language=Python]
init_p_list = [[[3,1]], [[1,1]]]
init_knots_list = [[[[0, 0, 0, 0, 1, 1, 1, 1], 
                     [0, 0, 1, 1]]], 
                   [[[0, 0, 1, 1], 
                     [0, 0, 1, 1]]]]
\end{lstlisting}
The orders and knot vectors for the model after order elevation are given by the following variables. 
\begin{lstlisting}[language=Python]
p_list_ele = [[[3, 3]], [[3, 3]]]
knots_list_ele = [[[[0, 0, 0, 0, 1, 1, 1, 1], 
                    [0, 0, 0, 0, 1, 1, 1, 1]]], 
                  [[[0, 0, 0, 0, 1, 1, 1, 1], 
                    [0, 0, 0, 0, 1, 1, 1, 1]]]]
\end{lstlisting}
Next, we pass the multilevel design information to the instance \lstinline{des2ana}. To keep the vertical patch straight in the axial direction, the horizontal and vertical coordinates are aligned along the axial direction using the \lstinline{set_cp_align} method.
\begin{lstlisting}[language=Python]
des2ana.set_init_knots_by_field(init_p_list, 
                                init_knots_list)
des2ana.set_order_elevation_by_field(p_list_ele, 
                                     knots_list_ele)
des2ana.set_knot_refinement()
des2ana.set_cp_align(field=0, align_dir_list=[1])
des2ana.set_cp_align(field=2, align_dir_list=[1])
\end{lstlisting}
Furthermore, we can create the shape optimization model with moving intersections by passing the \lstinline{nmopt} and \lstinline{des2ana} instances to the \lstinline{ShapeOptGroupMint} class and create the optimization problem.
\begin{lstlisting}[language=Python]
model = ShapeOptGroupMint(cpdesign2analysis=des2ana,
                          nonmatching_opt=nmopt)
model.init_parameters()
prob = om.Problem(model=model)
\end{lstlisting}
Then the optimizer is configured using \lstinline{prob.driver.options}. The optimization problem is set up with \lstinline{prob.setup()} and solved by \lstinline{prob.run_driver()}. Using the SLSQP optimizer with a tolerance of $10^{-9}$, the optimized geometry after 22 iterations is shown in Figure \ref{subfig:tbeam-geom-opt}. A cross-sectional view of the T-beam is illustrated in Figure \ref{subfig:tbeam-opt-slice}, indicating that the vertical patch moves to the center of the top patch, thereby minimizing the internal energy of the T-beam for the given load and boundary conditions with sufficiently small errors. The optimized configuration in Figure \ref{fig:tbeam-opt} also demonstrates that the T-junction between the top and vertical patches is well preserved.

\subsection{Tube under internal pressure} \label{subsec:example-tube}
\begin{figure}[!htbp]
    \centering
    \begin{subfigure}[!htbp]{0.49\textwidth}
        \includegraphics[width=\textwidth]{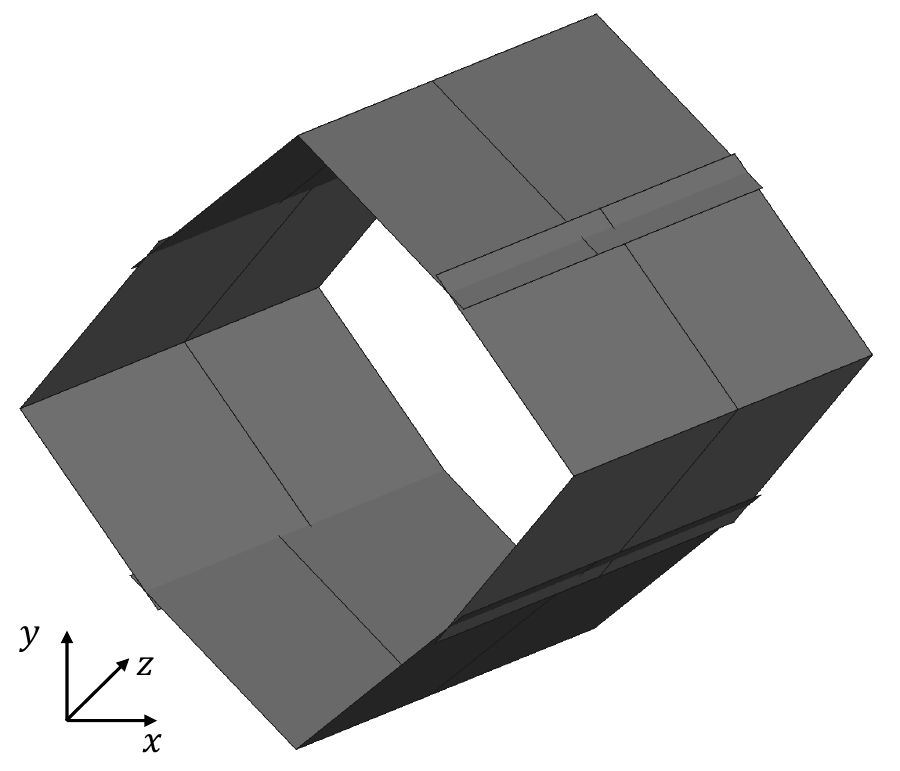}
        \caption{}
        \label{subfig:tube-geom-init}
    \end{subfigure}
    \hfill
    \begin{subfigure}[!htbp]{0.49\textwidth}
        \includegraphics[width=\textwidth]{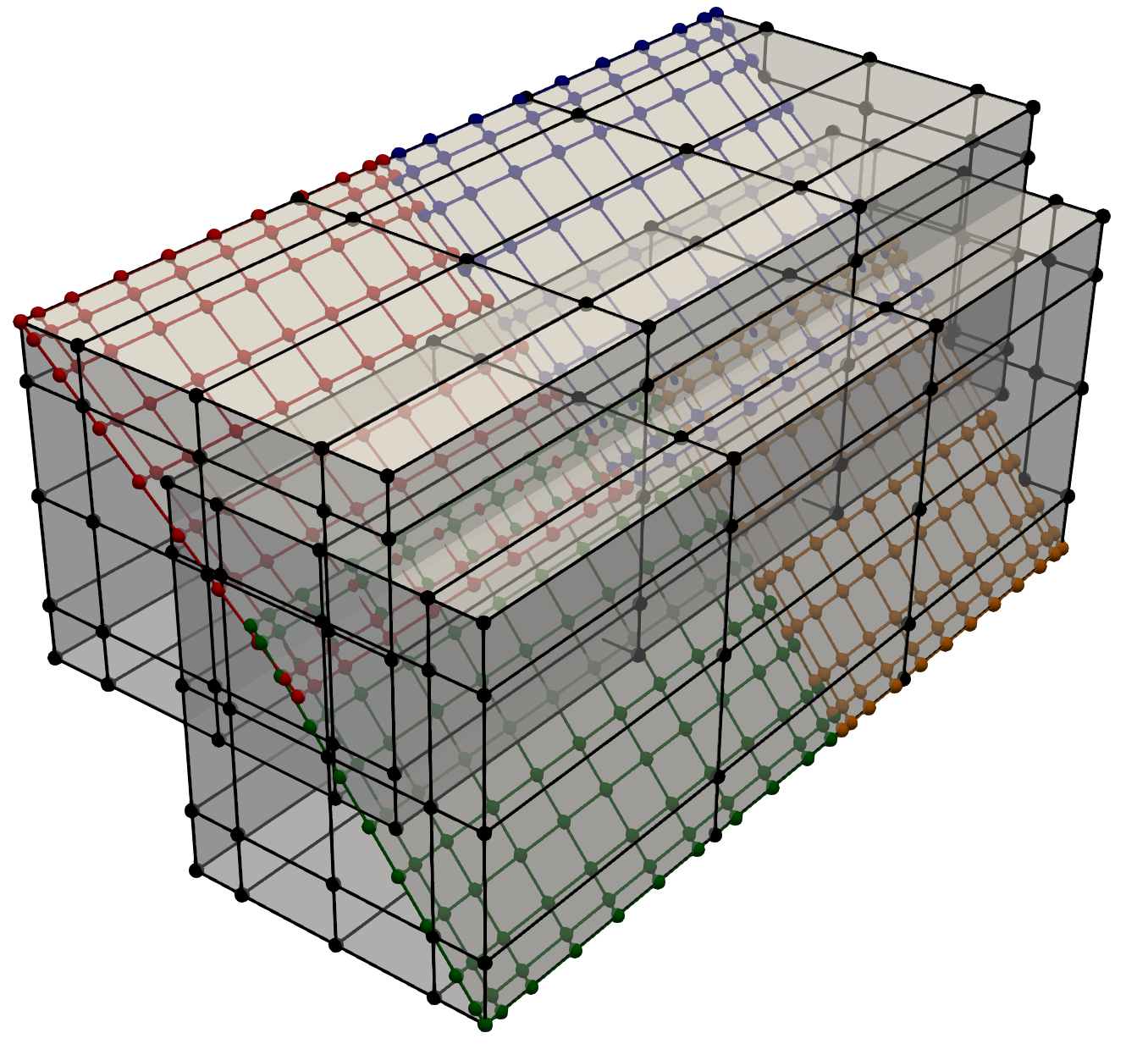}
        \caption{}
        \label{subfig:tube-geom-ffd}
    \end{subfigure}
    \caption{(a) The baseline design of a tube geometry, with a quarter of the tube modeled by four non-matching B-spline patches. (b) Two FFD blocks are employed, one for each pair of B-spline surfaces with an edge intersection. Relative movement is allowed between the two FFD blocks.}
    \label{fig:tube-geom}
\end{figure}

A tube under internal pressure is considered in this section. The initial design of the tube is shown in Figure \ref{subfig:tube-geom-init}. To validate the GOLDFISH framework, we model a quarter of the initial tube using four non-matching B-spline patches, with three differentiable intersections between the upper and lower pair of patches. Each pair of shell patches is embedded in one FFD block to maintain their edge intersection, while the intersections between the two FFD blocks are allowed to move during the optimization process. This setup enables the analytical optimal solution, which is a cylindrical tube. The configuration of the FFD blocks and their associated discretizations are illustrated in Figure \ref{subfig:tube-geom-ffd}. The optimization of this tube combines both the FFD-based approach and the moving intersections method.

\begin{figure}[!htbp]
    \centering
    \begin{subfigure}[!htbp]{0.49\textwidth}
        \includegraphics[width=\textwidth]{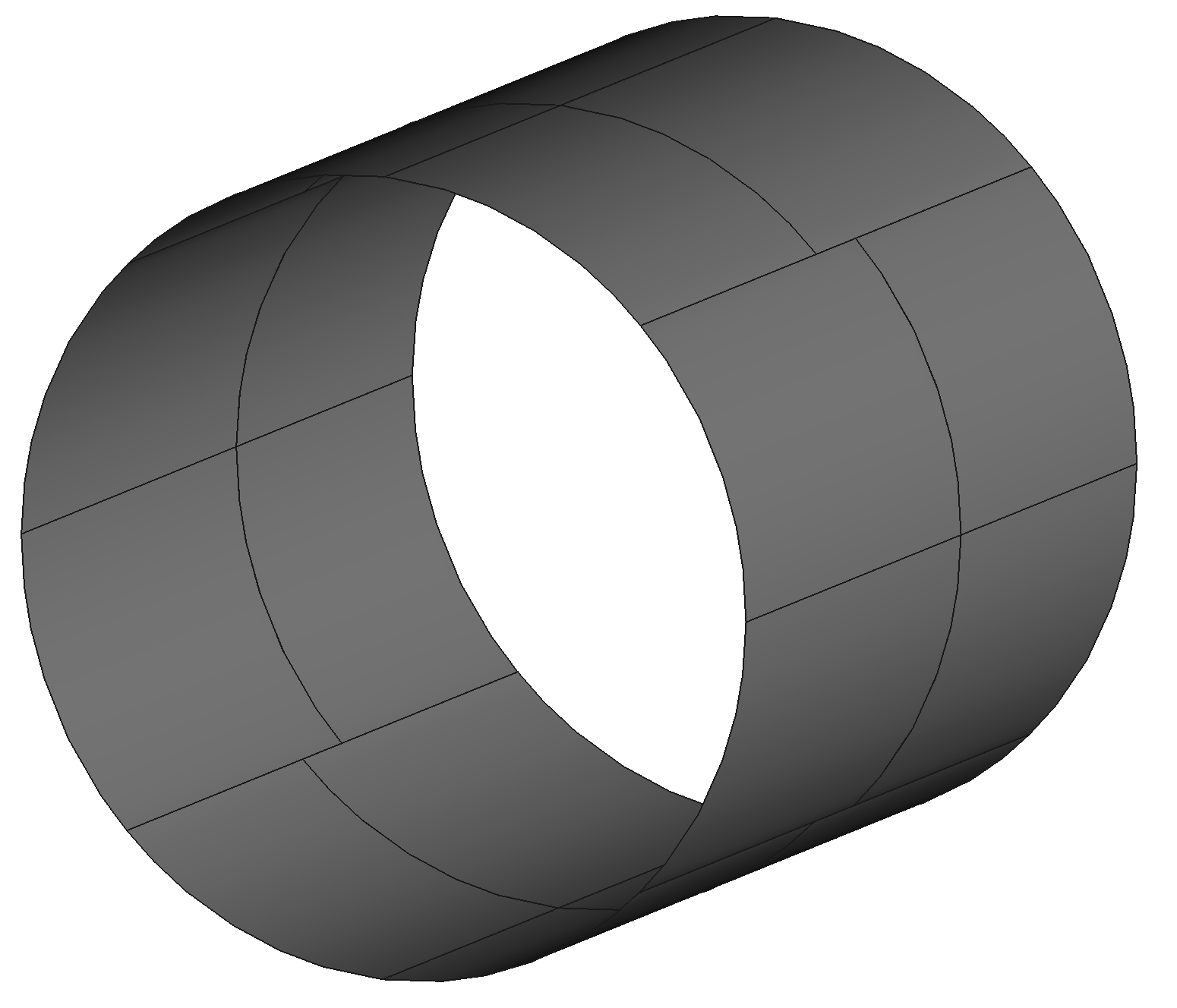}
        \caption{}
        \label{subfig:tube-geom-opt}
    \end{subfigure}
    \hfill
    \begin{subfigure}[!htbp]{0.49\textwidth}
        \includegraphics[width=\textwidth]{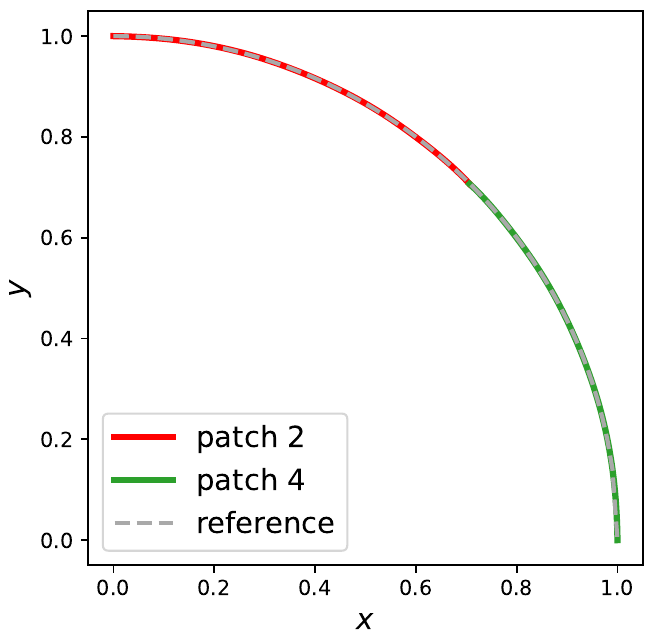}
        \caption{}
        \label{subfig:tube-slice}
    \end{subfigure}
    \caption{(a) The resulting geometry of the tube with minimum internal energy. Surface intersections between the spline patches transit to edge intersections in the optimized design. (b) Comparison of the optimized geometry with an exact cylindrical tube in the cross-sectional view.}
    \label{fig:tube-opt}
\end{figure}

The following code snippets illustrate the setup of the tube optimization problem using GOLDFISH. As in previous examples, we begin by creating instances of \lstinline{OCCPreprocessing} and \lstinline{NonmatchingOptFFD} for the imported CAD geometry. For both FFD blocks, we optimize the vertical and horizontal coordinates of their control points. Since the axial direction of the tube is the $z$ direction, the optimization field is set as [0,1] for both FFD blocks to update their $x$ and $y$ coordinates. Shell patches with indices 0 and 1 are embedded in the first FFD block, while the remaining two patches are embedded in the second FFD block. Thus, the \lstinline{opt_surf_inds} is defined as [[0,1], [2,3]]. This information is passed to the non-matching problem using the method \lstinline{set_shopt_surf_inds_multiFFD}.
\begin{lstlisting}[language=Python]
opt_field = [[0,1],[0,1]]
opt_surf_inds = [[0,1], [2,3]]
nm_opt.set_shopt_surf_inds_multiFFD(opt_field, opt_surf_inds)
nm_opt.set_geom_preprocessor(preproc)
\end{lstlisting}
We then create two lists containing the knot vectors \lstinline{shopt_ffd_knots_list} and control points \lstinline{shopt_ffd_control_list} for the trivariate B-spline solids used to define the FFD blocks, as demonstrated in Figure \ref{subfig:tube-geom-ffd}. The definitions of these lists are omitted for clarity. The method \lstinline{set_shopt_multiFFD} is used to obtain the data of the FFD blocks.
\begin{lstlisting}[language=Python]
nm_opt.set_shopt_multiFFD(shopt_ffd_knots_list,
                          shopt_ffd_control_list)
\end{lstlisting}                                
Next, the control points of both FFD blocks are aligned in the axial direction using the method \lstinline{set_shopt_align_CP_multiFFD}. The control points on the faces of the FFD blocks that lie on the symmetric planes are fixed using \lstinline{set_shopt_pin_CP_multiFFD}. Additionally, the method \lstinline{set_shopt_regu_CP_multiFFD} is employed to prevent self-penetration of the FFD blocks.
\begin{lstlisting}[language=Python]
nm_opt.set_shopt_align_CP_multiFFD(ffd_ind=0, 
                                   align_dir=[[2],[2]])
nm_opt.set_shopt_align_CP_multiFFD(ffd_ind=1, 
                                   align_dir=[[2],[2]])
nm_opt.set_shopt_pin_CP_multiFFD(ffd_ind=0, pin_dir0=[0,0], 
                                 pin_side0=[[0],[0]])
nm_opt.set_shopt_pin_CP_multiFFD(ffd_ind=1, pin_dir0=[1,1], 
                                 pin_side0=[[0],[0]])
nm_opt.set_shopt_regu_CP_multiFFD()
\end{lstlisting}
The remainder of the optimization setup is identical to the previous examples, which involves creating quadrature meshes for the intersections and defining PDE residuals for the Kirchhoff--Love shells. The OpenMDAO optimization problem is then created using \lstinline{om.Problem}. The SNOPT optimizer is employed for this problem with a tolerance of $10^{-2}$. The converged tube geometry after 142 iterations is shown in Figure \ref{subfig:tube-geom-opt}. A cross-sectional view of a quarter of the tube is displayed in Figure \ref{subfig:tube-slice} and compared with an exact quarter circle for validation. The surface intersections between the two pairs of shell patches move to the edges of the patches, forming a cylindrical tube to achieve the optimal shape. This demonstrates that both the FFD-based and moving intersections approaches can work simultaneously in shape optimization for non-matching shell structures.
% \begin{lstlisting}[language=Python]
% nm_opt.create_mortar_meshes(preproc.mortar_nels)
% nm_opt.mortar_meshes_setup(preproc.mapping_list, 
%        preproc.intersections_para_coords, 
%        penalty_coefficient, transfer_mat_deriv=2)
% preproc.check_intersections_type()
% preproc.get_diff_intersections()
% nm_opt.create_diff_intersections(num_edge_pts=None)
% \end{lstlisting}

\subsection{Application to aircraft wings} \label{subsec:example-wings}
This section presents the optimization results for aircraft wings to demonstrate the applicability of GOLDFISH to complex shell structures.

\subsubsection{Wing shape optimization using FFD-based approach} \label{subsubsec:example-evtol-wing-ffd}
\begin{figure}[!htbp]
    \centering
    \begin{subfigure}[!htbp]{0.49\textwidth}
        \includegraphics[width=\textwidth]{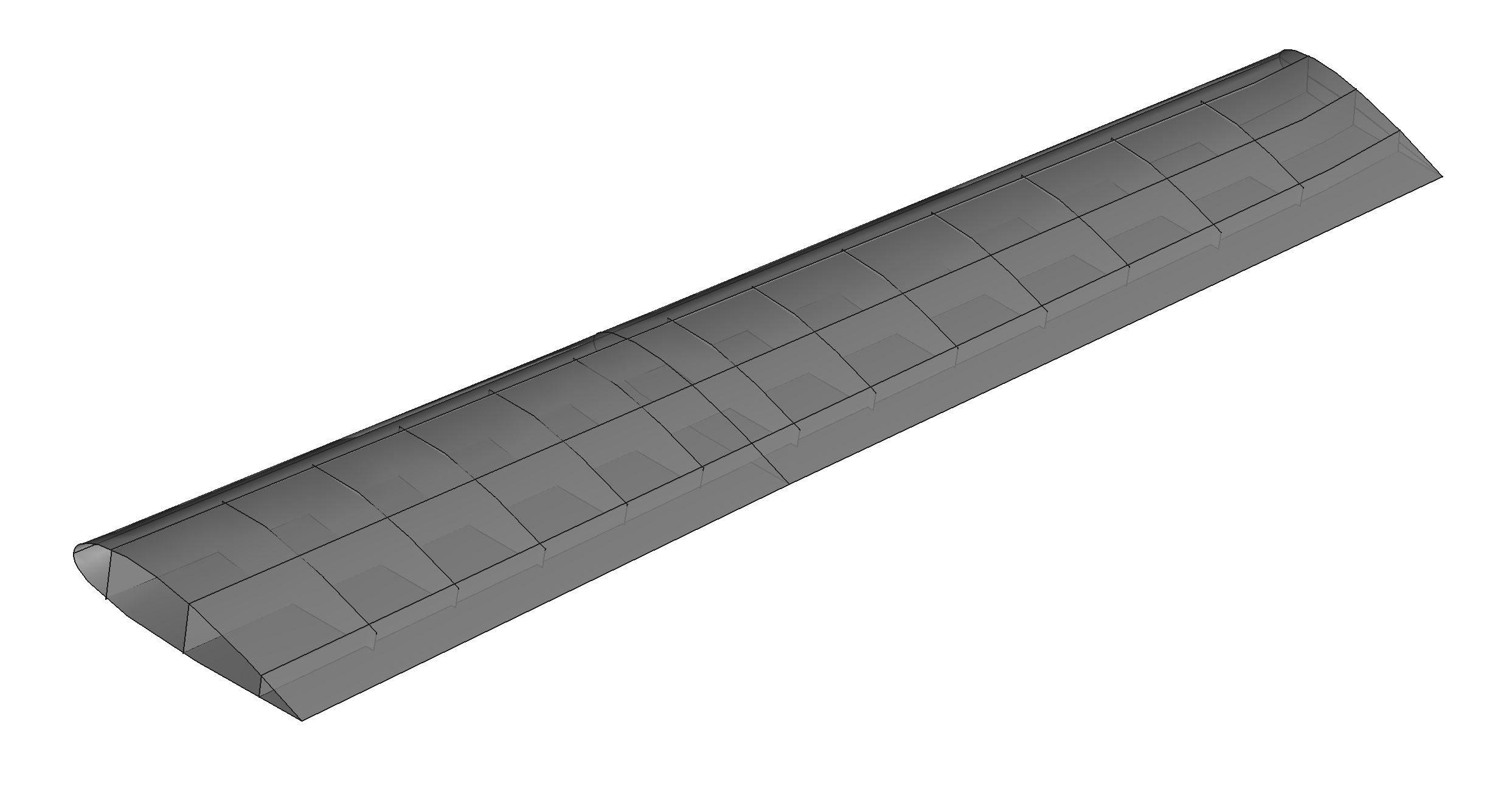}
        \caption{}
        \label{subfig:evtol-wing-geom-init-ffd}
    \end{subfigure}
    \hfill
    \begin{subfigure}[!htbp]{0.49\textwidth}
        \includegraphics[width=\textwidth]{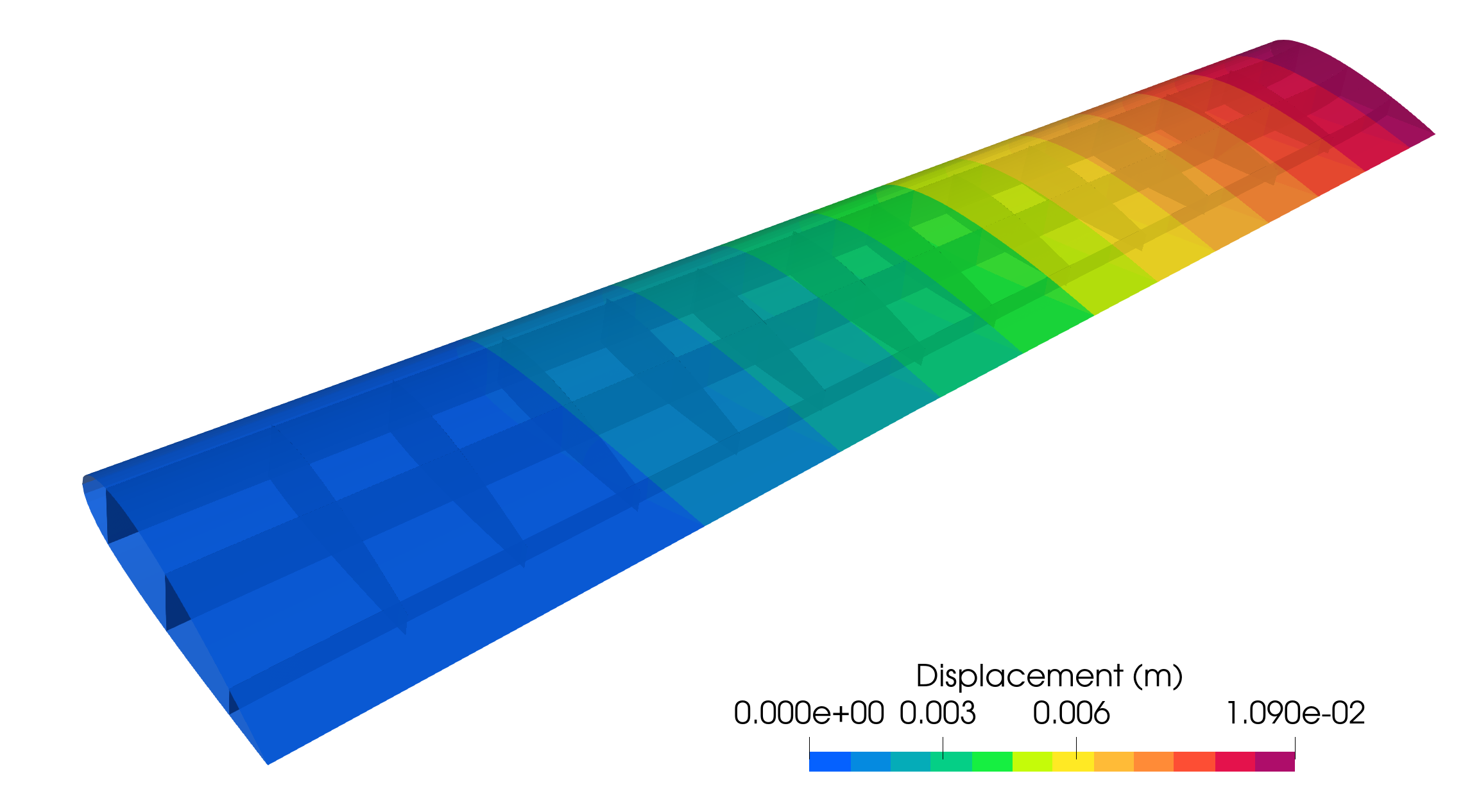}
        \caption{}
        \label{subfig:evtol-wing-disp-init-ffd}
    \end{subfigure}
    \caption{(a) The initial CAD geometry of an eVTOL wing. (b) The displacement field of the eVTOL wing under a distributed load.}
    \label{fig:evtol-wing-init-ffd}
\end{figure}

In the first application, we use the FFD-based approach to optimize the shape of an electric vertical take-off and landing (eVTOL) aircraft wing. The CAD geometry of the eVTOL wing, shown in Figure \ref{subfig:evtol-wing-geom-init-ffd}, is composed of 21 B-spline patches and 87 intersections. The wing is clamped at the root and subjected to an upward-facing distributed pressure. The material properties used are those of aluminum, with Young's modulus of 68 GPa and Poission's ratio of 0.35. All shell patches have a thickness of 3 mm. The contour plot of the wing displacement in the baseline design is shown in Figure \ref{subfig:evtol-wing-disp-init-ffd}.

\begin{figure}[!htbp]
    \centering
    \begin{subfigure}[!htbp]{0.49\textwidth}
        \includegraphics[width=\textwidth]{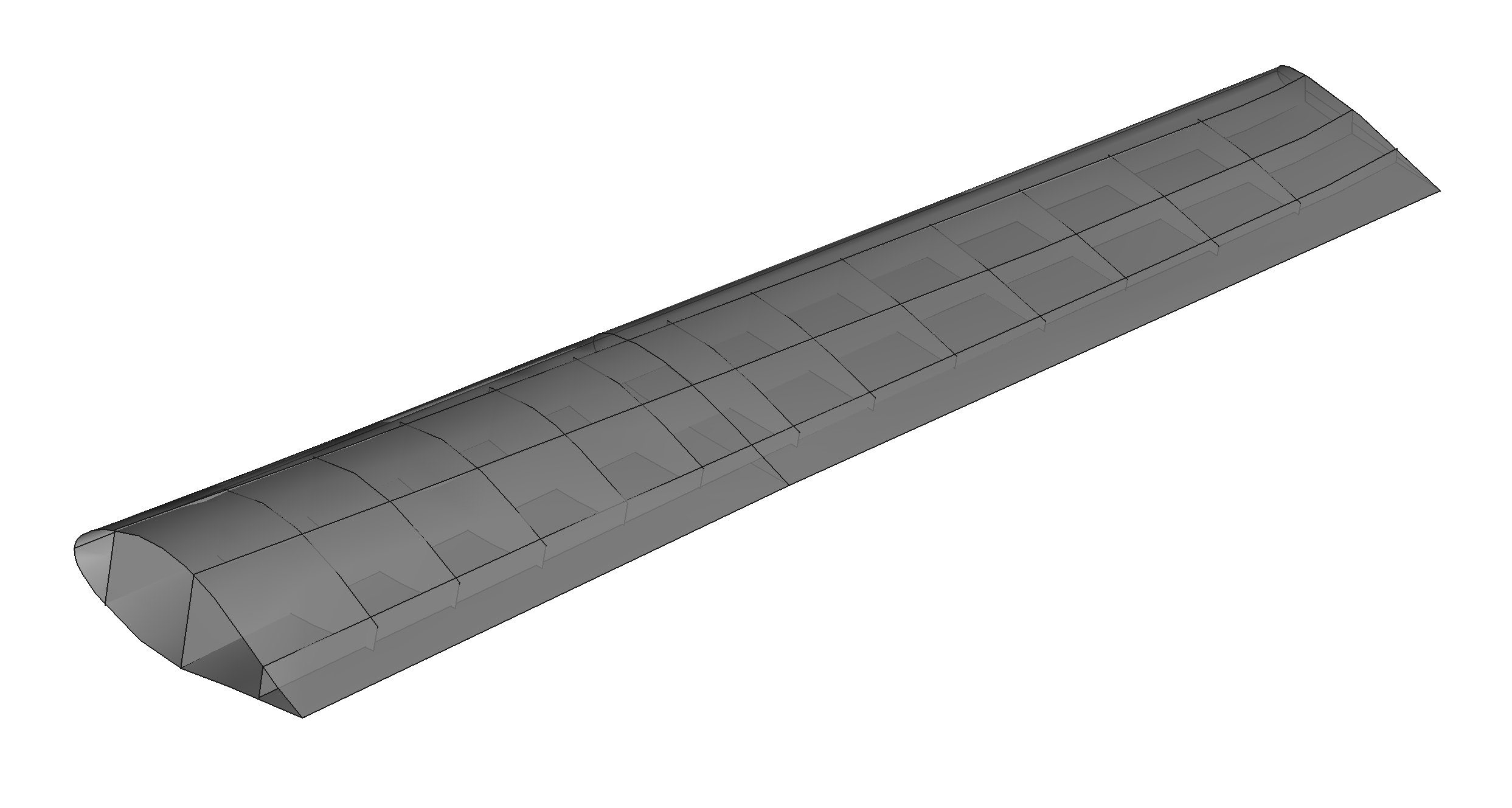}
        \caption{}
        \label{subfig:evtol-wing-geom-opt-ffd}
    \end{subfigure}
    \hfill
    \begin{subfigure}[!htbp]{0.49\textwidth}
        \includegraphics[width=\textwidth]{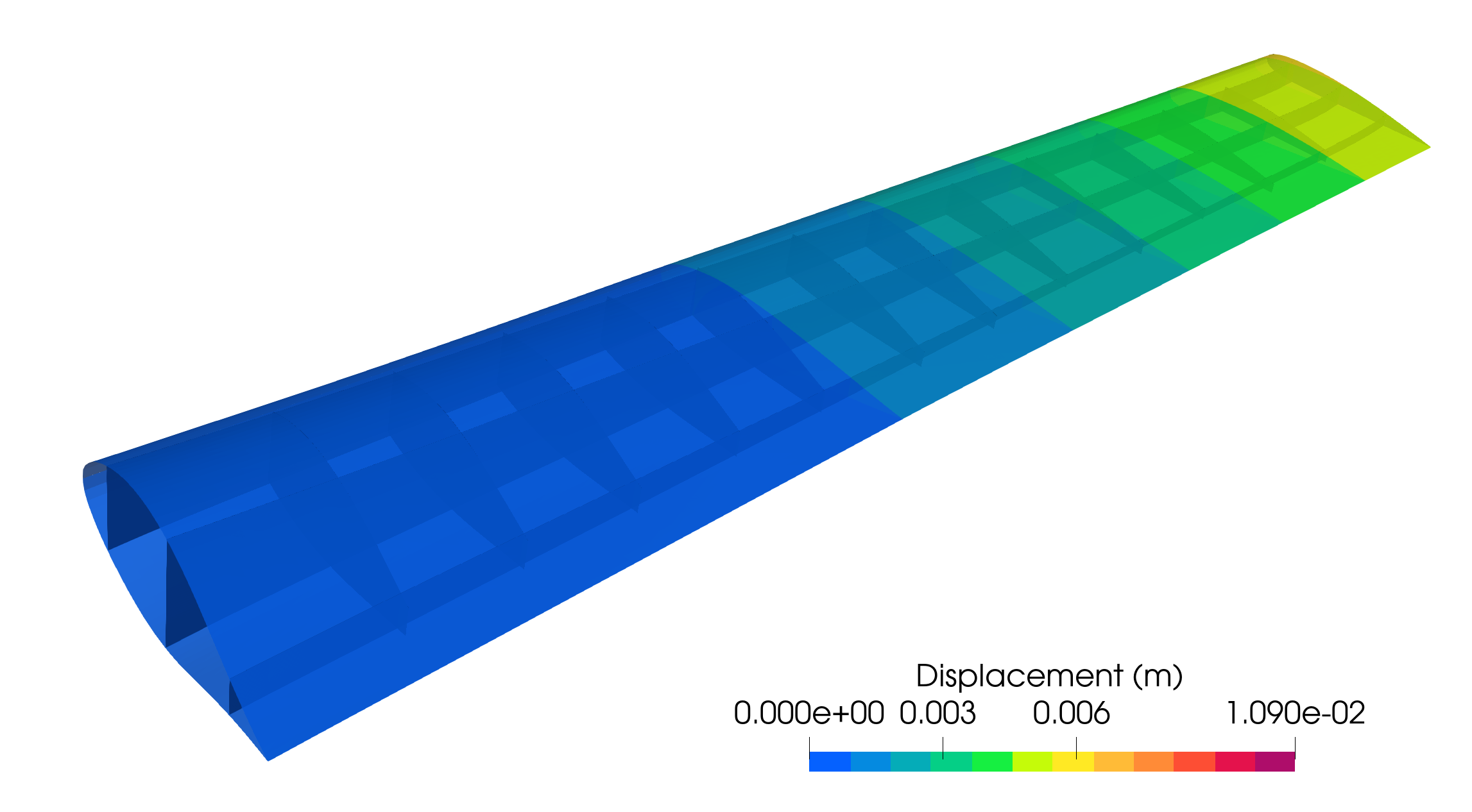}
        \caption{}
        \label{subfig:evtol-wing-disp-opt-ffd}
    \end{subfigure}
    \caption{(a) The optimized CAD geometry of the eVTOL wing using the FFD-based approach. (b) The displacement field of the optimized eVTOL wing.}
    \label{fig:evtol-wing-opt-ffd}
\end{figure}

To perform shape optimization for the eVTOL wing, the entire geometry is embedded in an FFD B-spline block for shape update. The control points of the FFD block in the vertical direction serve as the design variables. The internal energy of the wing is minimized under a constant volume constraint. A regularization term is added to the objective function to smooth the gradient changes of the surface in the vertical direction, preventing oscillatory shapes. The optimized geometry is shown in Figure \ref{subfig:evtol-wing-geom-opt-ffd}. The cross-section of the wing root becomes wider to provide greater support under the distributed pressure, while the wing tip narrows in the vertical direction to compensate for the increased volume at the wing root. The corresponding displacement contour plot of the optimized wing is displayed in Figure \ref{subfig:evtol-wing-disp-opt-ffd}, where the maximum displacement magnitude is noticeably decreased. The internal energy of the optimized wing is reduced by 49.1\% compared to the baseline design. 

\subsubsection{Internal structures shape optimization} \label{subsubsec:example-evtol-wing-mint}
\begin{figure}[!htbp]
    \centering
    \begin{subfigure}[!htbp]{0.49\textwidth}
        \includegraphics[width=\textwidth]{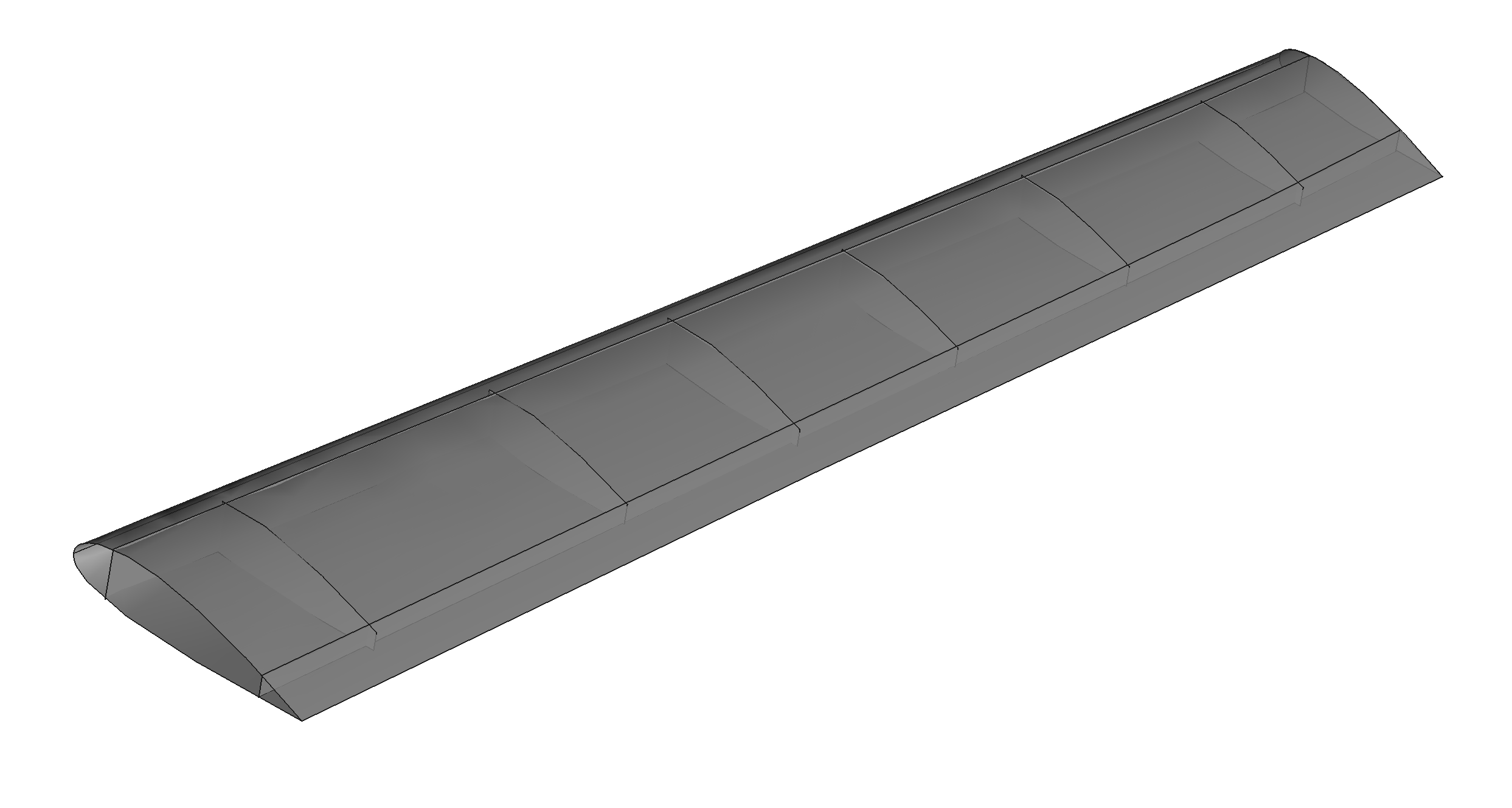}
        \caption{}
        \label{subfig:evtol-wing-geom-init-mint}
    \end{subfigure}
    \hfill
    \begin{subfigure}[!htbp]{0.49\textwidth}
        \includegraphics[width=\textwidth]{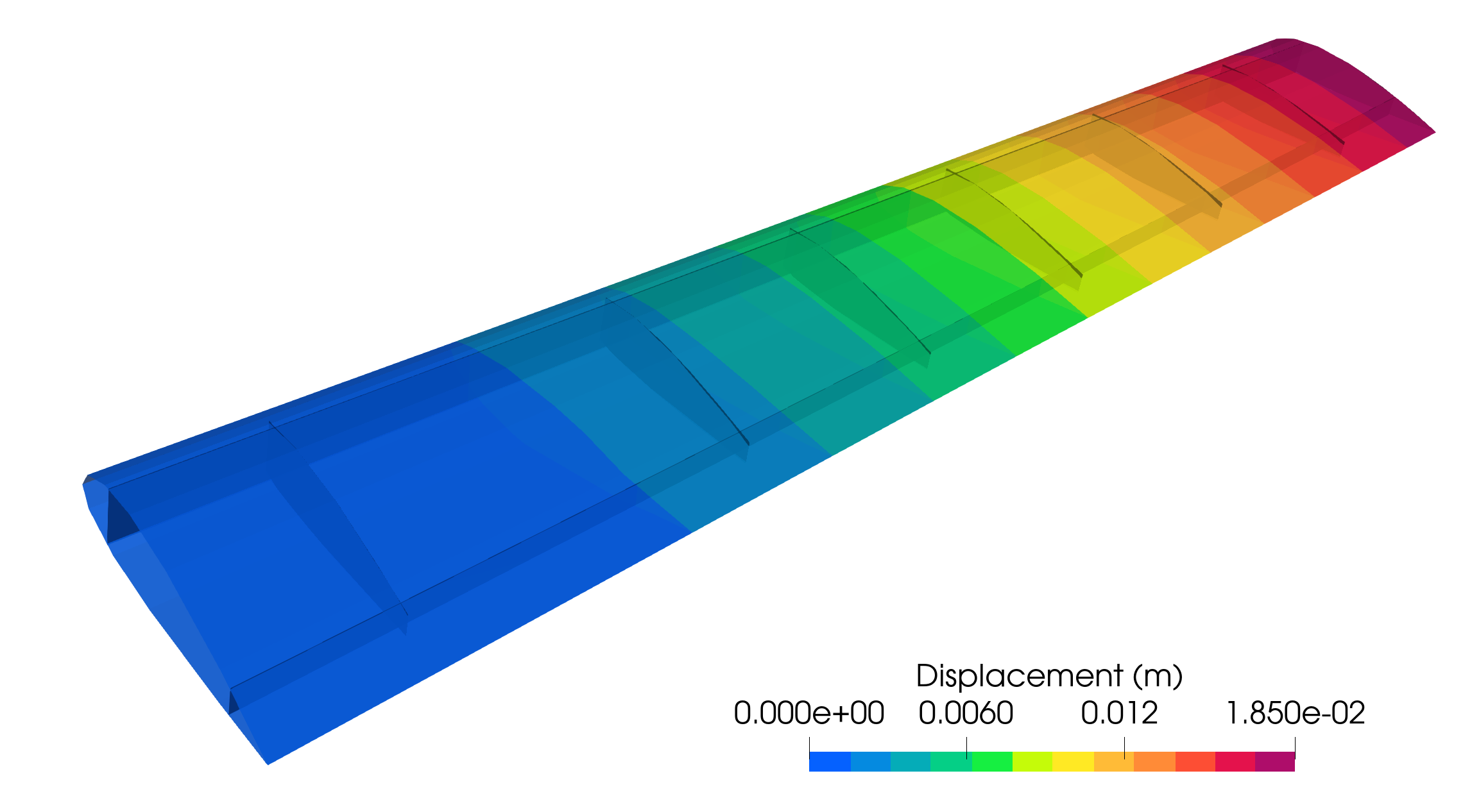}
        \caption{}
        \label{subfig:evtol-wing-disp-init-mint}
    \end{subfigure}
    \caption{(a) The baseline design of an eVTOL wing. (b) The displacement of the initial eVTOL wing.}
    \label{fig:evtol-wing-init-mint}
\end{figure}

This section focuses on the optimization of the layout of internal structures of the eVOTL wing, where intersections between the outer skins, spars, and ribs are allowed to move during the optimization process. The initial geometry of the eVTOL wing is shown in Figure \ref{subfig:evtol-wing-geom-init-mint}, consisting of 11 B-spline patches and 32 intersections. Among these intersections, four of them are located on the edges of intersecting outer surfaces and therefore remain fixed. The remaining 28 intersections are movable to optimize the layout of the internal sub-structures. The boundary conditions, loading conditions, material properties, and geometric parameters are the same as those in Section \ref{subsubsec:example-evtol-wing-ffd}, and the displacement contour plot is shown in Figure \ref{subfig:evtol-wing-disp-init-mint}.

\begin{figure}[!htbp]
    \centering
    \begin{subfigure}[!htbp]{0.49\textwidth}
        \includegraphics[width=\textwidth]{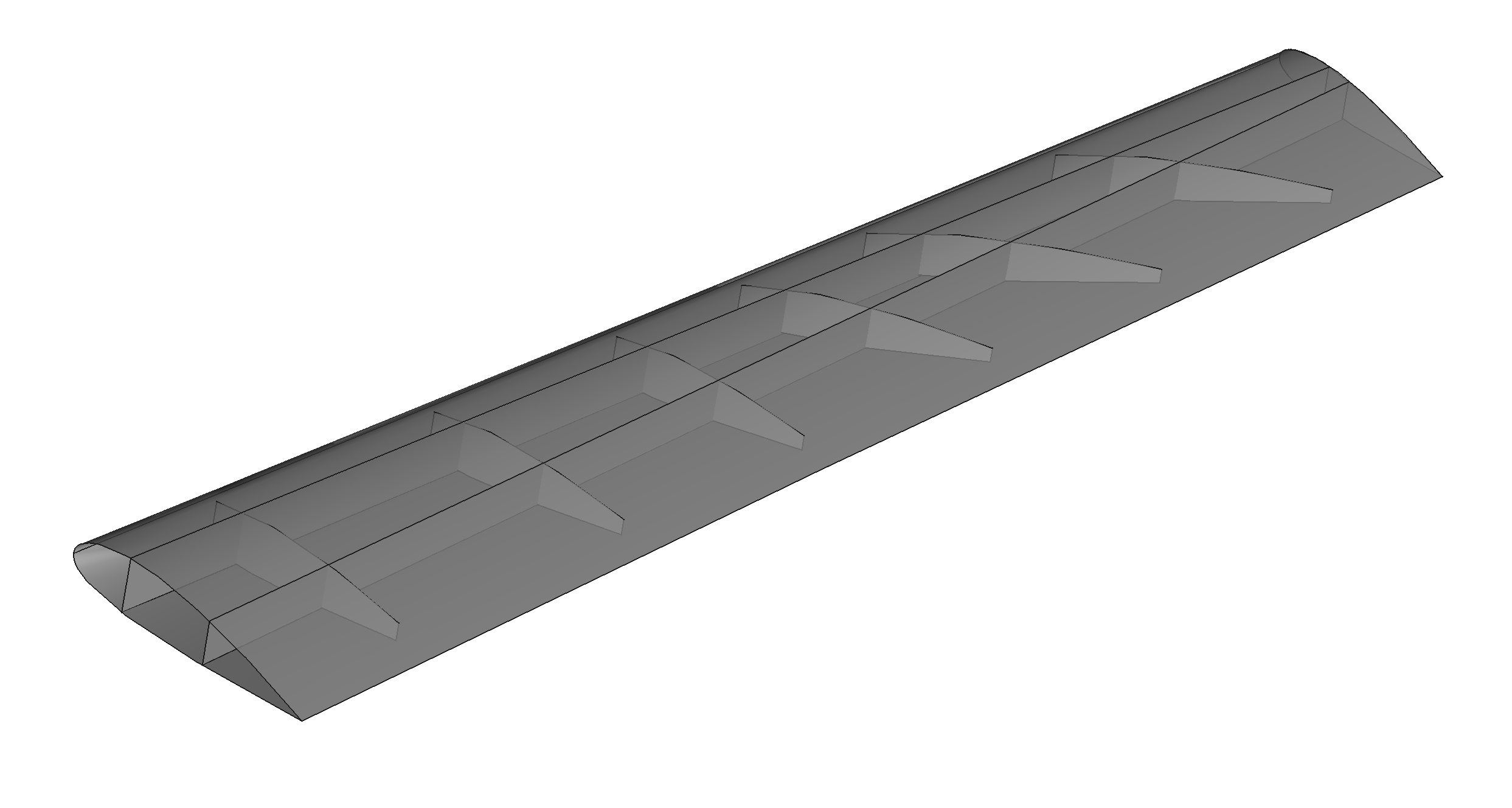}
        \caption{}
        \label{subfig:evtol-wing-geom-opt-mint}
    \end{subfigure}
    \hfill
    \begin{subfigure}[!htbp]{0.49\textwidth}
        \includegraphics[width=\textwidth]{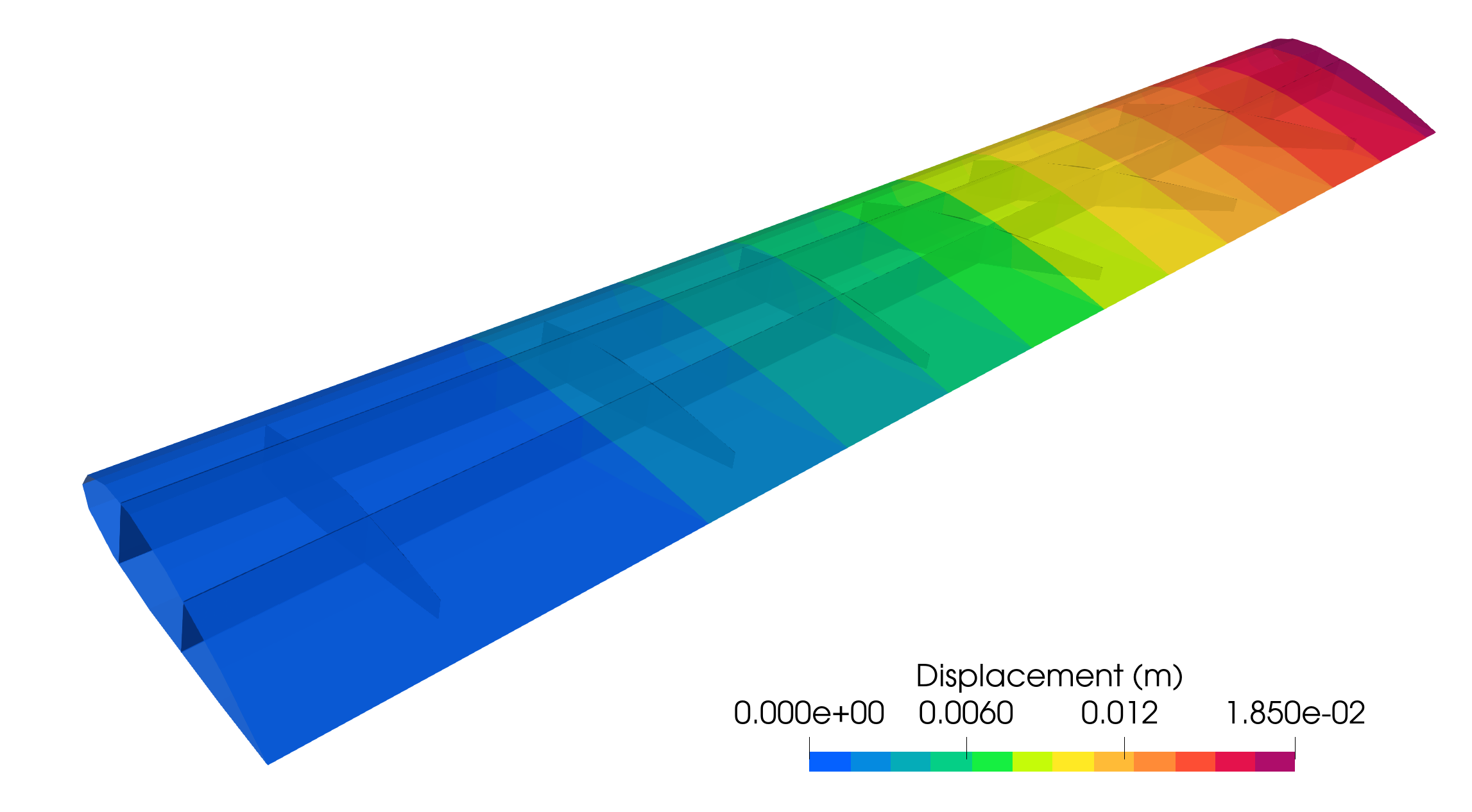}
        \caption{}
        \label{subfig:evtol-wing-disp-opt-mint}
    \end{subfigure}
    \caption{(a) The CAD geometry of the eVTOL wing with optimized internal ribs and spars. (b) The displacement field of the optimized eVTOL wing.}
    \label{fig:evtol-wing-opt-mint}
\end{figure}

In this application, the outer surfaces remain unchanged during the shape optimization. Internal spars are allowed to have rigid body translation within the envelope of the outer surfaces, while the ribs can translate and rotate along the lines from 15\% to 80\% of the distance from the leading edge to the trailing edge while remaining planar surfaces. Each rib has a volume constraint set to less than 1.5 times its initial value to prevent excessively stretched elements. The optimized geometry, which minimizes the internal energy, is shown in Figure \ref{subfig:evtol-wing-geom-opt-mint}, and the associated displacement field is visualized in Figure \ref{subfig:evtol-wing-disp-opt-mint}. The tip displacement is slightly decreased as the contour lines move toward the wing tip. The internal energy of the optimized wing is reduced by 5.2\% compared to the initial design. Additionally, the T-junctions between the internal structures and outer skins are maintained, despite the changes in the locations of internal structures.

The eVTOL wing shape optimization examples validate the applicability of GOLDFISH to complex shell structures with arbitrary intersections. The overall shape optimization exhibits a significant reduction in internal energy. Furthermore, the optimization of internal spars and ribs not only reduces the objective function but also leads to innovative layouts of internal structures. This section showcases two representative design scenarios for the eVTOL wing, further examples of wing design can be found in \cite[Section 6]{zhao2024automated} using the FFD-based approach and in \cite[Section 6]{zhao2024shape} with moving intersections.

\section{Conclusions} \label{sec:conclusions}
This paper introduces GOLDFISH, an open-source Python library for the shape optimization of complex shell structures. The structural analysis employs an isogeometric Kirchhoff--Love shell model coupled with a penalty formulation for patch intersections, and control points of shell patches are modified to update the shell shape during shape optimization. In this approach, FE mesh generation is no longer required in the optimization loop. This framework is developed based on FEniCS, leveraging its automatic differentiation and code generation capabilities to compute derivatives, thereby enabling gradient-based design optimization. This code framework accepts B-spline or NURBS-based CAD geometries as input and returns the optimized geometry in the same format, streamlining the workflow from geometry design through structural analysis to shape optimization.

The modular design of this framework inherited from OpenMDAO is presented along with the essential components for FFD-based shape optimization and the moving intersections approaches. A suite of benchmark problems, accompanied by code implementations, is provided to examine the effectiveness of the framework. Furthermore, applications to aerospace structures are presented, demonstrating innovative designs for the layout of internal spars and ribs in an eVTOL wing using this GOLDFISH open-source code. Nevertheless, conditions from other disciplines, such as aerodynamics \cite{van2023solver, van2024enforcing}, electric motors \cite{scotzniovsky2024geometric}, and noise control \cite{gill2023applicability}, need to be considered for practical designs. From an engineering perspective, shape optimization for aerospace structures with stress and buckling constraints represents a promising direction for framework development. While this GOLDFISH code offers basic IGA-based capabilities for multi-patch shape optimization of Kirchhoff-Love shell structures, additional enhancement of the implemented mathematical formulation and numerical algorithms can be performed under the current framework. For example, implementation of the membrane locking-free Kirchhoff--Love shell using mixed formulation, special quadrature schemes, B-bar type projection \cite{koschnick2005discrete, echter2013hierarchic, bouclier2013efficient, greco2018reconstructed, bieber2018variational, casquero2022removing, casquero2023overcoming, sauer2024simple, mathews2024computationally} can be considered. The B-bar type techniques can be implemented by constructing a projection operator to modify the strain-displacement matrix. Furthermore, special quadrature schemes are more involved within the FEniCS-based code environment, and mixed formulations would require substantial modifications to the optimization framework due to additional unknowns. These directions offer opportunities for future development by integrating various disciplines and techniques to enhance the framework for more practical shell structure designs. The code framework is maintained on the GitHub repository \cite{goldfish-code}.

\backmatter

%\bmhead{Supplementary information}

%If your article has accompanying supplementary file/s please state so here. 

%Authors reporting data from electrophoretic gels and blots should supply the full unprocessed scans for key as part of their Supplementary information. This may be requested by the editorial team/s if it is missing.

%Please refer to Journal-level guidance for any specific requirements.

% \bmhead{Acknowledgments}
% -- Acknowledgements

\section*{Acknowledgments}\label{sec:Ackonwledgements}
Han Zhao was supported by NASA grant number 80NSSC21M0070 when preparing the original submission. We thank Dr. David Kamensky at the University of California San Diego for helpful discussions on FEniCS implementation.

%\section*{Declarations}

%Some journals require declarations to be submitted in a standardised format. Please check the Instructions for Authors of the journal to which you are submitting to see if you need to complete this section. If yes, your manuscript must contain the following sections under the heading `Declarations':

%\begin{itemize}
%\item Funding
%\item Conflict of interest/Competing interests (check journal-specific guidelines for which heading to use)
%\item Ethics approval 
%\item Consent to participate
%\item Consent for publication
%\item Availability of data and materials
%\item Code availability 
%\item Authors' contributions
%\end{itemize}

%\noindent
%If any of the sections are not relevant to your manuscript, please include the heading and write `Not applicable' for that section. 

\bibliography{main}% common bib file
%% if required, the content of .bbl file can be included here once bbl is generated
%%\input sn-article.bbl

%% Default %%
%%\input sn-sample-bib.tex%

\end{document}